\documentclass[12pt, leqno]{article}
\usepackage{amssymb,amsmath}
\usepackage[pdftex]{color}
\usepackage{latexsym}
\usepackage[pdftex]{graphicx}
\newtheorem{theorem}{Theorem}
\newtheorem{proposition}{Proposition}
\newtheorem{definition}{Definition}
\newtheorem{corollary}{Corollary}
\newtheorem{remark}{Remark}
\newtheorem{lemma}{Lemma}
\newtheorem{example}{Example}
\newcommand{\tx}[1]{\mbox{\;{#1}\;}} 
\newcommand{\N}{\mathbb{N}}
\newcommand{\R}{\mathbb{R}^n}
\newcommand{\rz}{\mathbb{R}}

\newcommand{\Div}{\hbox{div}\:}
\newcommand{\p}{\partial}
\numberwithin{equation}{section}
\begin{document}
\title{Second Domain Variation for Problems with Robin Boundary Conditions}
\pagestyle{myheadings}
\maketitle
\centerline{\scshape Catherine Bandle}
\medskip
{\footnotesize
\centerline{Mathematische Institut, Universit\"at Basel,}
\centerline{Rheinsprung 21, CH-4051 Basel, Switzerland}
} 
\medskip
\centerline{\scshape Alfred Wagner}
\medskip
{\footnotesize
\centerline{Institut f\"ur Mathematik, RWTH Aachen  }
\centerline{Templergraben 55, D-52062 Aachen, Germany}}
\maketitle
\begin{abstract}
In this paper, the first and second domain variation for functionals related to elliptic boundary and eigenvalue problems with Robin boundary conditions is computed. Extremal properties of the ball among nearly spherical domains of given volume are derived. The discussion leads to a Steklov eigenvalue problem. As a by-product a general characterization of the optimal shapes is obtained.
\end{abstract}
%%%%%%%%%%%%%%%%%%%%%%%%%%%%%%%%%%%%%%%%%%%
\newline
\newline
\centerline {MSC2010 49K20, 49J20, 34L15, 35J20, 35N25.}
%%%%%%%%%%%%%%%%%%%%%%%%%%%%%%%%%%%%%%%%%%%
%%%%%%%%%%%%%%%%%%%%%%%%%%%%%%%%%%%%%%%%%%%
\section{Introduction}
\label{intro}
We consider elliptic boundary and eigenvalue problems with Robin boundary conditions. 
We are interested in the question how functionals, such as the associated Robin energy or the principal eigenvalue, depend on the domain. For this purpose, we shall use the method of domain variations, which - in the spirit of calculus - studies the functionals under an infinitesimal change of the domain. 
\newline
\newline
This technique has a long history and can be traced back to Hadamard \cite{Ha}. 
Fifty years later Garabedian and Schiffer \cite{GaSc53} computed the first and second domain
variation for several functionals such as the first eigenvalue of the Dirichlet-Laplace operator, 
the virtual mass and the Green's function. They chose special perturbations and obtained 
convexity theorems. 
Subsequent to the work of Garabedian and Schiffer, D. Joseph \cite{Jo67} computed higher variations of the eigenvalues and studied the behavior of the spectrum under shear and stretching. Grinfeld \cite{Gr10} computed the eigenvalues of a polygon. For a long time this topic has rather been neglected. In the last years it has attracted considerable interest. New developments and new applications are found in the inspiring books by Henry \cite{He}, Khludnev and Sokolowski \cite{KhSo97} and Pierre and Henrot \cite{HePi05}, where further references are given. 
\newline
\newline
Classical domain functionals, like the first eigenvalue of the Dirichlet-Laplace operator or the torsional rigidity, satisfy isoperimetric inequalities in the class of domains of prescribed volume. They are known as the Rayleigh-Faber-Krahn and the St. Ve\-nant - P\'olya inequality \cite{PoSz51}. We are interested in finding similar inequalities. In contrast to the classical domain functionals, the Robin energy contains an additional boundary energy term.
To our knowledge, the effect of this boundary integral on the second domain variation hasn't been explored yet.
\newline
\newline
We shall compute the first and second variations of the domain functionals by means of a change of variable to keep the domain fixed. Then we study the perturbed problem on this fixed domain. 
\newline
\newline
The first variation is a simple. It provides a necessary condition for extremal domains in terms of an overdetermined elliptic problem. The first eigenvalue of the ball is a candidate for an extremal domain. The same is true for other energies if the solutions of the Euler Lagrange equation are radial. This is in accordance with the Bossel-Daners inequality \cite{Da06}, which states that among all domains of given volume the ball yields a local minimum of the first eigenvalue, and by recent results by Bucur and Giacomini \cite{BuGia13} for nonlinear problems. As a by-product, we obtain a local monotonicity property which improves slightly the one in \cite{Gi}.
\newline
\newline
The second variation was studied by many authors, see for instance \cite{KhSo97} and Pierre and Novruzi \cite{NoPi02} and references cited in their work. The abstract structure seems to be well understood. However, the strict positivity (coercivity), which is necessary for the minimality property of a domain, remained in general a challenging open problem. 
\newline
Fujii \cite{Fu90} showed that for the torsion problem in the plane, the second variation of the energy is positive for area preserving perturbations. In \cite{KaKu14} H. Kasumba and K. Kunisch consider an optimal control problem, which is related to our energy. They minimize a prescribed cost functional subject to a Poisson equation with Dirichlet data. They compute the second domain variation.
There are other approaches to estimate optimality. Instead of using the second domain variation, one can use the Fraenkel asymmetry function. For the Dirichlet eigenvalue problem for the Laplace operator this was investigated in \cite{BrPr}.
\newline
\newline
The aim of this paper is to discuss the sign of the second variation for volume preserving perturbations. The two main ingredients are a device by Simon \cite{Si80} and a Steklov eigenvalue problem.
In contrast to problems with Dirichlet boundary conditions, the second variation of the surface plays a crucial role, see also \cite{DaMi20}. First attempts to compute the first and second variations can be found in \cite{BaWa10} and \cite{BaWa15}.
\newline
\newline
Our paper is organized as follows.
\newline
\newline
First we present, for the reader's convenience, some concepts and tools from Riemannian geometry, which will be needed for handling the shape derivatives of the boundary integrals. We then discuss useful properties of the vector fields which are related to volume preserving perturbations. In Section 3 we describe in full detail the energies and the Rayleigh quotients of the perturbed problems, expressed in the original domain. The first variations are derived in Section 4. They lead to a characterization of an optimal domain. In Section 5 the shape derivatives of the solutions in the perturbed domains will be discussed. They will play an essential role for the sign of the second variation. Section 6 is devoted to the lengthly computations of the second variation. Applications to problems in nearly spherical domains of fixed volume are investigated in Section 7. As a surprise we find that the sign of the second variation for the ball depends on the sign of the surface energy. We compare our approach with Garabedian and Schiffer's formula for the second variation of the principal eigenvalue of the Laplacian with Dirichlet boundary conditions. We show that the ball is a local minimum. For the sake of completeness we give the formula for the second variation of the energy in case of Dirichlet boundary conditions.
%%%%%%%%%%%%%%%%%%%%%%%%%%%%%%%%%%%%%%%%
%%%%%%%%%%%%%%%%%%%%%%%%%%%%%%%%%%%%%%%%
\section{Preliminaries}
\label{prel}
%%%%%%%%%%%%%%%%%%%%%%%%%%%%%%%%%%%%%%%%
%%%%%%%%%%%%%%%%%%%%%%%%%%%%%%%%%%%%%%%%
\subsection{Geometry of  surfaces}\label{curvature}
%%%%%%%%%%%%%%%%%%%%%%%%%%%%%%%%%%%%%%%%
In this section we recall some ideas about surfaces in $\R$.
Throughout this paper we will use the following notation. Let $\Omega\subset\R$ be a  bounded 
$C^{2,\alpha}$-domain, and let  $x:=(x_1,x_2,\dots,x_n)$ denote a  point in $\R$. 
The expression $x\cdot y$ stands for the Euclidean scalar product of two vectors $x$ and 
$y$ in $\R$ and $|x|=(x\cdot x)^{1/2}$.
\newline
At every point $P\in \partial \Omega$ there exists a neighborhood $\cal{U}_P$ and a Cartesian coordinate system with an orthonormal basis $\{e^{i}\}_{i=1}^n$ centered at  $P$, such that $e^n$ points in the direction of the outer normal $\nu$ and $e^i $, $i=1,\dots n-1$ lie in the tangent space of $P$. The coordinates with respect to this basis will be denoted by $(\xi_1,\xi_2,\dots,\xi_n)$. Moreover, we assume that 
\begin{eqnarray*}
\Omega\cap {\cal{U}}_P= \{ \xi\in {\cal{U}}_P: \xi_n<F(\xi_1,\xi_2,\dots,\xi_{n-1})\}, \quad F\in C^{2,\alpha}.
\end{eqnarray*}
With  this choice of coordinates clearly $F(0)=0$
and $F_{\xi_i}(0)=0$ for $i=1,2,\dots,n-1$.
For short we set $\xi':=(\xi_1,\xi_2,\cdots,\xi_{n-1})$ which ranges in $\mathcal{U}':=\mathcal{U}_P\cap\{\xi_n=0\} $.
\newline
\newline
In ${\cal{U}}_P \cap \p \Omega$ the boundary is represented by $x(\xi') =(\xi_1,\xi_2,\dots, \xi_{n-1},F(\xi'))$ and the unit  outer normal $\tilde\nu(\xi')=(\tilde\nu_1,\tilde\nu_2.\dots,\tilde\nu_n)$ with respect to the $\xi'-$coordinate system, is given by
\begin{eqnarray*}
\tilde\nu(\xi')= \frac{(-F_{\xi_1}, -F_{\xi_2},\dots,-F_{\xi_{n-1}},1)}{\sqrt{1+|\nabla' F|^2}},
\end{eqnarray*}
where $\nabla'$ stands  for the gradient in $\rz^{n-1}$. 
In this paper we shall use the Einstein convention, where repeated indices are understood to be summed from $1$ to $n-1$ or from $1$ to $n$, respectively. The vectors $x_{\xi_i}$, $i=1,2,\dots,n-1$ span the tangent space. The metric tensor of $\p\Omega$ is denoted by $g_{ij}$ and its inverse by $g^{ij}$. We have
\begin{eqnarray*}
g_{ij} =x_{\xi_i}\cdot x_{\xi_j}&=&\delta_{ij} + F_{\xi_i}F_{\xi_j}\qquad\hbox{and}\\ 
g^{ij}&=& \delta_{ij} -\frac{F_{\xi_i}F_{\xi_j}}{1+|\nabla' F|^2}.
\end{eqnarray*}
The surface element of $\p\Omega$ is 
\begin{eqnarray*}
dS=\sqrt{{\rm{det} }g_{ij}}d\xi'=\sqrt{1+|\nabla' F |^2}d\xi'.
\end{eqnarray*} 
Observe that any vector $v$  can be represented in the form 
\begin{eqnarray}\label{vrepresentation}
v= g^{ij}(v\cdot x_{\xi_i})x_{\xi_j} +(v\cdot \nu)\nu, \: i,j=1,2,\dots, n-1.
\end{eqnarray}
Let $f\in C^1(\cal{U}_P)$ and let $\tilde f(\xi'):=f(\xi)|_{\p \Omega}=f(\xi',F(\xi'))$. The {\sl tangential gradient} of  $f$ at a boundary point is defined as
\begin{eqnarray}\label{tangrad}
\nabla ^{\tau} \tilde f := g^{ij}\frac{\p \tilde f}{\p \xi_j}x_{\xi_i}.
\end{eqnarray}
Let us write for short $\p_i:= \frac{\p}{\p \xi_i}$. Then
\begin{eqnarray}\label{tander}
\p^*_i\tilde{f}:=g^{ij}\p _j \tilde f
\end{eqnarray}
is the tangential derivative on $\p\Omega$.
For a smooth vector field $ v:\p \Omega \to \R$, which is not necessarily tangent to $\p \Omega$, we define the {\sl tangential divergence} by
\begin{eqnarray}\label{tandivmet}
\Div_{\p\Omega}v:= g^{ij}\tilde v_{\xi_i}\cdot x_{\xi_j}, \tx{where} \tilde v=v(\xi', F(\xi')).
\end{eqnarray}
By \eqref{tander} this can also be written as
\begin{eqnarray}\label{tandivmet1}
\Div_{\p\Omega}v=\partial_{j}^{*}\tilde v\cdot x_{\xi_{j}}.
\end{eqnarray}
For $i=1,\dots,n-1$ the numbers $\kappa_i$ are the principal curvatures of $\p\Omega$ in the point $P$. We define
\begin{eqnarray*}
H:=\frac{1}{n-1}\sum_{i=1}^{n-1}\kappa_i
\end{eqnarray*}
This is the {\sl mean curvature of $\partial\Omega$} at $P$. For an arbitrary point $(\xi', F(\xi'))$ on $\p \Omega$ it is given by
\begin{eqnarray*}
H(\xi')= (n-1)^{-1}\frac{\p}{\p\xi_i}\left(-\frac{F_{\xi_i}(\xi')}{\sqrt{1+|\nabla' F|^2}}\right).
\end{eqnarray*}
Observe that
\begin{align}\label{meanc}
\Div_{\p \Omega} \nu =(n-1)H.
\end{align}
In particular we have $H=\frac{1}{R}$ if $\partial \Omega=\partial B_R$, where $B_R$ denotes the ball of radius $R$ centered at the origin.
\newline
\newline
Another way of defining the tangential gradient and the tangential divergence is by projection onto the tangent space of $\partial \Omega$. Let $x\in \partial\Omega$ and let $T_{x}\partial\Omega$ be the tangent space of $\partial\Omega$ in $x$. Then we define
\begin{eqnarray*}
P\: :\: \R\to T_{x}\partial\Omega\qquad v\to P(v):=v-(v\cdot\nu)\nu.
\end{eqnarray*}
We will also write $v^{\tau}:=P(v)$. From \eqref{vrepresentation} the gradient $\nabla f$ in $\R$ can be decomposed as
\begin{align*}
\nabla f=g^{ij}(\nabla f\cdot x_{\xi_i})x_{\xi_j} + (\nabla f\cdot \nu)\nu.
\end{align*}
Note that
$
\nabla f\cdot x_{\xi_i}= \p_if + \p_nfF_{\xi_i} = \p_i\tilde f.
$
Hence
\begin{align}\label{tgradient}
\nabla^\tau \tilde f =\nabla f -(\nabla f\cdot \nu)\nu=P(\nabla f).
\end{align}
As in \cite{Giu84} some computations will be shorter if we introduce  the $i$-th component of the tangential gradient 
\begin{align*}
\delta_if= \p_if -\nu_i\p_jf \nu_j.
\end{align*} 
At the origin we have $\delta_i f=\p^*_i f= \p_i f$. In general 
$\delta_i f$ and $\p^*_if$ are different, more precisely
\begin{align}\label{umrech}
\delta_k = (x_{\xi_j}\cdot e^k) \p^*_j=(\p^*_jx\cdot e^k)\p_j.
\end{align}
\newline
In the same way we show that for any smooth vector field $v:\overline{\Omega}\to\R$
\begin{eqnarray}\label{tandiv}
\Div_{\partial\Omega}v=\Div v-\nu\cdot D_v \nu:=\partial_{i}v_{i}-\nu_{j}\partial_{j}v_{i}\nu_{i}=\delta_j \tilde v_j.
\end{eqnarray}
At the origin we have $\Div_{\p\Omega}v =\p_i v_i=\p_i^*\tilde v_i$, $i=1,\dots n-1$.
\newline
\newline
We will frequently use integration by parts on $\partial\Omega$. Let $f\in C^1(\partial\Omega)$ and $v\in C^{0,1}(\partial\Omega,\R)$ then the Gauss theorem on surfaces holds:
\begin{eqnarray}\label{tanpi}
\oint_{\p\Omega}f\:\Div_{\partial\Omega}v\:dS
=-\oint_{\p\Omega}v\cdot\nabla^\tau f\:dS
+(n-1)\oint_{\p\Omega}f(v\cdot\nu)\:H\:dS.
\end{eqnarray}
This formula can also be written in the form
\begin{eqnarray}\label{tanpi2}
\oint_{\p\Omega}f\:\delta_{j}v_{j}\:dS
=-\oint_{\p\Omega}v_{j}\delta_{j}f\:dS
+(n-1)\oint_{\p\Omega}f(v\cdot\nu)\:H\:dS.
\end{eqnarray}
%%%%%%%%%%%%%%%%%%%%%%%%%%%%%%%%%%%%%%%%%%%%%%%%%%
\subsection{Domain perturbations}
\label{DP}
\subsubsection{Volume element}
\label{DPVE}
%%%%%%%%%%%%%%%%%%%%%%%%%%%%%%%%%%%%%%%%%%%%%%%%%%%
From now on we consider a bounded domain $\Omega\subset \R$ in the H\" older class $C^{2, \alpha}$. We consider a family of perturbations
$(\Omega_t)_{t}$ of the form
\begin{align}\label{yvar}
\Omega_t:=\{y=x+tv(x)+\frac{t^2}{2}w(x)+o(t^2):x \in \Omega,\: t\tx{small}\},
\end{align}
where $v=(v_1(x),v_2(x),\dots,v_n(x))$ and $w=(w_1(x),w_2(x),\dots,w_n(x))$ are $C^{2, \alpha}$  vector fields. With $o(t^2)$ we denote all terms such that $\frac{o(t^2)}{t^2}\to 0$ as $t\to 0$. This regularity assumption justifies the computations carried out in the sequel. 
\newline
\newline
The  Jacobian matrix corresponding to the transformation  $y(t,\Omega):=y$ in \eqref{yvar}
is - up to second order terms - given by
\begin{align*}
I + t D_v +\frac{t^2}{2}D_w, \tx{where} (D_v)_{ij}=\partial_j v_i \tx{and} \partial_j=\partial/ \partial x_j.
\end{align*}
By Jacobi's formula we have for small $t$
\begin{eqnarray}\label{jform}
J(t)&:=&\tx{det}(I+tD_v+\frac{t^2}{2}D_w)\\\nonumber&=&1+t\:\Div v+\frac{t^2}{2}\left((\Div v)^2-D_v : D_v+\Div\: w\right)+o(t^2).
\end{eqnarray}
Here we have used the notation
\begin{eqnarray*}
D_v : D_v:=\partial_{i}v_j\partial_{j}v_{i}.
\end{eqnarray*}
Thus, $y(t,\Omega)$ is a diffeomorphism for $t\in ]-t_0,t_0[$ and $t_0$ sufficiently small.
\newline
\newline
\noindent  {\sl Throughout this paper we shall consider diffeomorphisms $y(t,\Omega)$ as described above.}
\newline
\newline
Later on we will be interested in volume preserving transformations. From 
\begin{eqnarray*}
|\Omega_t| &=& \int_\Omega J(t)\:dx=|\Omega| +t\int_\Omega \tx{div} v\:dx \\ 
\qquad&&+ \frac{t^2}{2}\int_\Omega((\Div v)^2-D_v : D_v+\Div\:w)\:dx+o(t^2) 
\end{eqnarray*}
it follows, that $y(t,\Omega)$ is {\sl  volume preserving of the first order} if
\begin{align}\label{volume1}
\int_\Omega \Div v\:dx =0.
\end{align}
It is {\sl volume preserving of the second order}, if in addition to \eqref{volume1} it satisfies
\begin{align}\label{volume2}
\int_\Omega((\Div v)^2-D_v : D_v+\Div \:w)\:dx=0.
\end{align}
For volume preserving transformations of the second order we have
%%%%%%%%%%%%%%%%%%%%%%%%%%%%%%%%%%
\begin{lemma}\label{volpres1} Let $v\in C^{0,1}(\Omega,\R)$. 
Then
\begin{eqnarray*}
\int_{\Omega}\left((\Div v)^2-D_v : D_v+\Div\:w\right)\:dx=0
\end{eqnarray*}
is equivalent to
\begin{align}\label{volcomp1}
\oint_{\p\Omega}(v\cdot\nu) \Div v\:dS-\oint_{\p\Omega}v_{i}\:\partial_{i}v_j\nu_j\:dS+\oint_{\p\Omega}(w\cdot \nu)\:dS=0.
\end{align}
\end{lemma}
%%%%%%%%%%%%%%%%%%%%%%%%%%%%%%%%%llllllllllllllllllllll
{\bf{Proof}}
Integration by parts gives
\begin{eqnarray*}
\int_{\Omega}D_v : D_v\:dx
&=&
-
\int_\Omega v_{j}\:\partial_{j}\left(\Div v\right)\:dx
+
\oint_{\p\Omega}v_{j}\:\partial_{j}v_{i}\:\nu_{i}\:dS\\
&=&
\int_{\Omega}\left(\Div v\right)^2\:dx
-
\oint_{\p\Omega}(v\cdot\nu)\:\Div v\:dS
+
\oint_{\p\Omega}v_i \:\partial_iv_j\:\nu_j\:dS.
\end{eqnarray*}
This proves \eqref{volcomp1}. 
\hfill{$\rlap{$\sqcap$}\sqcup$}
\newline\newline
%%%%%%%%%%%%%%%%%%%%%%%%%%%%%%%%%llllllllllllllllllllll
\begin{remark}
The presence of $w$ is crucial, because otherwise the class of perturbations is too limited. For instance consider $B_1 \subset \mathbb{R}^2$ and let $\Omega_t$ 
be a rotated ball. Thus,
\begin{equation*}
y=
\left(
\begin{matrix}
\cos(t) & -\sin(t) \\
\sin(t)  & \cos(t)
\end{matrix}
\right) x.
\end{equation*}
Then (for small $t$) $y= x +t (-x_2,x_1) -\frac{t^2}{2}(-x_1,x_2)+o(t^2).$ It is easy to see that
the first order approximation $x+t(-x_2,x_1)$ is not volume preserving of the second order.
\end{remark}
%%%%%%%%%%%%%%%%%%%%%%%%%%%%%%%%%%%%%%
\begin{remark}\label{remarkw}
For any given $v$ we can always find a vector $w$ such that \eqref{volcomp1} is satisfied.
\end{remark}
%%%%%%%%%%%%%%%%%%%%%%%%%%%%%%%%%%%%%%%%%
\subsubsection{Surface element}
\label{subsubse}
%%%%%%%%%%%%%%%%%%%%%%%%%%%%%%%%%%
In this part we shall compute the surface element of  $\p \Omega_t$. Let $x(\xi')$, 
$\xi'\in \mathcal{U}'_P$ be local coordinates of $\p \Omega$ introduced in 
Section \ref{curvature}. Then $\p \Omega_t$ is represented locally by 
\begin{align*}
\{y(\xi'):=x(\xi')+t\tilde v(\xi')+\frac{t^2}{2}\tilde{w}(\xi'):\xi'\in \mathcal{U}'_P\},
\end{align*}
As before we write $\tilde v(\xi') := v(\xi',F(\xi'))$ and $\tilde w(\xi') := w(\xi',F(\xi'))$.
We set
\begin{eqnarray*}
g_{ij}&:=&x_{\xi_i}\cdot x_{\xi_j}, \\ 
a_{ij}&:=&x_{\xi_i}\cdot 
\tilde v_{\xi_j}+x_{\xi_j}\cdot \tilde v_{\xi_i},\\ 
b_{ij}&:=&2\tilde v_{\xi_i}\cdot \tilde v_{\xi_j}+\tilde w_{\xi_i}\cdot x_{\xi_j}+\tilde w_{\xi_j}\cdot x_{\xi_i}.
\end{eqnarray*}
We get
\begin{eqnarray*}
|dy|^2&:=&(g_{ij}+ta_{ij}+\frac{t^2}{2}b_{ij})d\xi_id\xi_j =:g_{ij}^td\xi_id\xi_j.
\end{eqnarray*}
Write for short $G:=(g_{ij})$, $G^{-1}:=(g^{ij})$,  $A:=(a_{ij})$, $B:=(b_{ij})$ and correspondingly $G^t:=(g^t_{ij}).$ Then the surface element on $\partial\Omega_t$ is
\begin{eqnarray*}
dS_y= \left(\mbox{det}G^t\right)^{1/2}d\xi'.
\end{eqnarray*}
Clearly
\begin{eqnarray*}
\sqrt{\mbox{det}G^t}= \sqrt{\mbox{det}G}\{\underbrace{\rm{det}(I+tG^{-1}A+\frac{t^2}{2}G^{-1}B)}_{k(x,t)}\}^{1/2}.
\end{eqnarray*}
Set 
\begin{eqnarray*}
\sigma_A:=\rm{trace}\:G^{-1}A, \: \sigma_B=\rm{trace}\:G^{-1}B \tx{and} \sigma_{A^2}= \rm{trace}\:(G^{-1}A)^2.
\end{eqnarray*}
For small $t$ the Taylor expansion yields
\begin{align*}
k(x,t)= 1+t\sigma_A
+\frac{t^2}{2}\left(\sigma_B +\sigma_A^2-\sigma_{A^2}\right)+o(t^2)
\end{align*}
and
\begin{align}\label{Wurzelk}
\sqrt{k(x,t)}= 1+\frac{t}{2}\sigma_A+\frac{t^2}{2}\left(\frac{1}{2}\sigma_B+\frac{1}{4}\sigma_A^2 -\frac{1}{2}
\sigma_{A^2} \right)+o(t^2).
\end{align}
In the sequel we shall use the notation 
\begin{eqnarray*}
m(t):= 1+\frac{t}{2}\sigma_A +\frac{t^2}{2}\left(\frac{1}{2}\sigma_B-\frac{1}{2}\sigma_{A^2}+\frac{1}{4}\sigma^2_A\right)+o(t^2).
\end{eqnarray*} 
Then the surface element of $\p \Omega_t$ reads as
\begin{eqnarray}\label{surface} 
dS_t = m(t)dS, \tx{where} dS \tx{is the surface element of $\p\Omega$}.
\end{eqnarray}
Our next goal is to find more explicit forms for the  expressions in $m(t)$.
It follows immediately from Section 2.1 that
\begin{align*}
\sigma_A=2g^{ij}\tilde v_{\xi_j}\cdot x_{\xi_i} =2\Div_{\p\Omega} v.
\end{align*}
The expression $\sigma_A$ has a geometrical interpretation. We find after a straightforward computation that
\begin{align}\label{SpurA}
\frac{1}{2}\sigma_A = \Div_{\p \Omega} v^{\tau} +(n-1)H( v\cdot \nu),
\end{align}
where $v^\tau$ is the projection of $v$ into the tangent space.
\newline
\newline
Moreover, a straightforward calculation leads to
\begin{align*}
\sigma_{A^2}&=g^{is}g^{kl}a_{sk}a_{li}=2(\p^*_i \tilde v\cdot x_{\xi_k})(\p^*_k\tilde v\cdot x_{\xi_i})+2(\p^*_i \tilde v\cdot x_{\xi_k})(\tilde v_{\xi_i}\cdot \p^*_k x),\\
\sigma_B&= 2g^{ij}\tilde v_{\xi_i}\cdot  \tilde v_{\xi_j} + 2\Div_{\p \Omega} w\\
&= 2(\p^*_s\tilde v\cdot x_{\xi_m})(\tilde v_{\xi_s}\cdot \p^*_mx) +2g^{is}(\tilde v_{\xi_i}\cdot\tilde  \nu)(\tilde v_{\xi_s}\cdot \tilde \nu) +2\Div_{\p\Omega} w.
\end{align*}
In the last expression we have used for $\tilde v_{\xi_k}$ the representation  \eqref{vrepresentation}.
Consequently
\begin{eqnarray}\label{ddmzero}
\ddot{m}(0)&=&\frac{1}{2}\sigma_B-\frac{1}{2}\sigma_{A^2}+\frac{1}{4}\sigma^2_A\\\nonumber&=& g^{is}(\tilde v_{\xi_i}\cdot \tilde \nu)(\tilde v_{\xi_s}\cdot \tilde \nu)+\Div_{\p \Omega} w- (\p^*_i \tilde v\cdot x_{\xi_k})(\p^*_k\tilde v\cdot x_{\xi_i})+  (\Div_{\p \Omega} v)^2.
\end{eqnarray}
Denote the surface area of $\Omega_t$ by
$\mathcal{S}_\Omega(t)$. Clearly its second variation is
\begin{eqnarray*}
\frac{d^2}{dt^2}S_{\Omega}(t)\vert_{t=0}=:\ddot{\mathcal{S}}_\Omega(0)=\oint_{\partial\Omega}\ddot{m}(0)\:dS.
\end{eqnarray*}
Then \eqref{ddmzero} together with \eqref{tanpi} implies, (see also \cite{Giu84,Li12})
\begin{eqnarray}\label{isodef}
\oint_{\partial\Omega}\ddot{m}(0)\:dS&=
\oint_{\partial\Omega}(\partial^{*}_{s}\tilde v\cdot\tilde \nu)\:(\partial_{s}\tilde v\cdot \tilde \nu)
-(\partial^{*}_{i}\tilde v\cdot x_{\xi_k})\:(\partial^{*}_{k}\tilde v\cdot x_{\xi_i})\\
\nonumber &
+
\left(\partial^{*}_{i}\tilde v\cdot x_{\xi_i}\right)^2\:dS
+(n-1)\int_{\p \Omega}( \tilde w\cdot \nu)H\:dS.
\end{eqnarray}
%%%%%%%%%%%%%%%%%%%%%%%%%%%%%%%%%%%%%%
\subsection{Computations for the ball}
\label{compball}
%%%%%%%%%%%%%%%%%%%%%%%%%%%%%%%%%%%%%%
In this subsection we simplify \eqref{isodef} for the special case $\partial\Omega=\partial B_R$. 
For simplicity we move the Cartesian coordinate system $\{e^{i}\}_{i=1}^n$ into the center of the ball. This transformation does not affect formula \eqref{isodef}. Note that in the radial case
\begin{eqnarray}\label{radder}
\nu= \frac{x}{R}\qquad\hbox{and}\qquad\delta_{i}\nu_{j}=\frac{1}{R}\left(\delta_{ij}-\nu_{i}\nu_{j}\right).
\end{eqnarray}
We now start with the evaluation of the different terms in \eqref{isodef}. We set $N:=(v\cdot \nu)$. We get
\begin{eqnarray*}
\ell_1:=(\partial^{*}_{s}\tilde v\cdot\tilde \nu)\:(\partial_{s}\tilde v\cdot \tilde \nu)&=\p^*_sN\p_sN-2(\p^*_s\tilde v\cdot \tilde \nu)(\tilde v\cdot \p_s \tilde \nu) +(\tilde v\cdot \p^*_s \tilde \nu)(\tilde v\cdot \p_s \tilde \nu)\\
&= |\nabla^\tau N|^2-2(\p^*_s\tilde v\cdot \tilde \nu)(\tilde v\cdot \p_s \tilde \nu) -R^{-2}|v^\tau|^2.
\end{eqnarray*}
By \eqref{radder} and \eqref{umrech}
\begin{eqnarray*}
-2(\p_s^*\tilde v\cdot \tilde \nu)(\tilde v\cdot \p_s \tilde \nu)&=&-\frac{2}{R} \p_s\tilde v_k\tilde \nu_k\tilde v_m(\p^*_s x\cdot e^m)\\
&=& -\frac{2}{R}[v_m\delta_mN -v_kv_m\p_s\nu_k(\p^*_s x\cdot e^m)]\\
&=& -\frac{2}{R}[(v\cdot \nabla^\tau N)- R^{-1}|v^\tau|^2].
\end{eqnarray*}
This, together with the Gauss theorem on surfaces \eqref{tanpi} implies
\begin{eqnarray}\label{intermed0}
&&\oint_{\partial B_R}\underbrace{(\p^*_s\tilde v\cdot\tilde \nu)\:(\p_s\tilde v\cdot \tilde \nu)}_{\ell_1}\:dS\\
\nonumber\qquad&=&\oint_{\partial B_R}(|\nabla^\tau N|^2 +\frac{2}{R}N\:\Div_{B_R}v
-\frac{2(n-1)}{R^2} N^2 + \frac{1}{R^2}|v^\tau|^2 )\:dS.
\end{eqnarray}
It is convenient to eliminate the last term in \eqref{intermed0}. If we replace in the Gauss formula  \eqref{tanpi}) $f$ by $N$ and use \eqref{radder}, we obtain
\begin{eqnarray*}
\frac{1}{R^2}\oint_{\partial B_R}|v^\tau|^2 \:dS&=&-\frac{1}{R}\oint_{\partial B_R}v_{j}\:\delta_{j}v_{i}\:\nu_{i}\:dS-\frac{1}{R}\oint_{\partial B_R}\Div_{\partial B_R}v\:N\:dS\\
&&+\frac{n-1}{R}\oint_{\partial B_R}N^2\:dS.
\end{eqnarray*}
%%%%%%%%%%%%%%%%%%%%%%%%%%%%%%%%%%%%%%%%%
With this remark we rewrite \eqref{intermed0} in the form
\begin{eqnarray}\label{intermed1}
\oint_{\partial B_R}\underbrace{(\p^*_s\tilde v\cdot\tilde \nu)\:(\p_s\tilde v\cdot \tilde \nu)}_{\ell_1}\:dS
&=&
\oint_{\partial B_R}\left(|\nabla^\tau N|^2 -\frac{(n-1)}{R^2} N^2\right)\:dS\\
\nonumber&&+
\frac{1}{R}\oint_{\partial B_R}\left(N\:\Div_{B_R}v -v_{j}\:\delta_{j}v_{i}\:\nu_{i}\right)\:dS.
\end{eqnarray}
Next we treat the second term. Observe that 
\begin{align*}
\ell_2=-(\partial^{*}_{i}\tilde v\cdot x_{\xi_k})\:(\partial^{*}_{k}\tilde v\cdot x_{\xi_i})=-\delta_s v_j \delta_j v_s.
\end{align*}
By \eqref{tanpi2} we find
\begin{eqnarray*}
\oint_{\partial B_R}\ell_2\:dS
=
\oint_{\partial B_R}v_j\:\delta_{i}\delta_{j}v_i\:dS
-
\frac{n-1}{R}\oint_{\partial B_R}v_j\:\delta_{j}v_i\nu_{i}\:dS.
\end{eqnarray*}
At this point it is important to note that $\delta_{i}\delta_{j}\ne\delta_{j}\delta_{i}$. In \cite{Giu84} (Lemma 10.7) the following relation is proved:
\begin{eqnarray*}
\delta_{i}\delta_{j}=\delta_{j}\delta_{i}+(\nu_{i}\delta_{j}\nu_k-\nu_{j}\delta_{i}\nu_k)\delta_{k}.
\end{eqnarray*}
For the ball this gives
\begin{align*}
(\nu_{i}\delta_{j}\nu_k-\nu_{j}\delta_{i}\nu_k)\delta_{k}=R^{-1}(\nu_i\delta_j-\nu_j\delta_i).
\end{align*}
Hence
\begin{eqnarray}\label{intermdiate3}
\oint_{\partial B_R}\ell_2\:dS
&=&
\oint_{\partial B_R}v_j\:\delta_{j}\delta_{i}v_i\:dS
-
\frac{n-2}{R}\oint_{\partial B_R}v_j\:\delta_{j}v_i\nu_{i}\:dS\\
\nonumber&&-
\frac{1}{R}\oint_{\partial B_R}N\:\Div_{\partial B_R}v\:dS.
\end{eqnarray}
We apply \eqref{tanpi2} to the first integral on the right side of \eqref{intermdiate3} and obtain
\begin{eqnarray*}
\oint_{\partial B_R}v_j\:\delta_{j}\delta_{i}v_i\:dS
=
-
\oint_{\partial B_R}(\delta_{i}v_{i})^2\:dS
+
\frac{n-1}{R}\oint_{\partial B_R}N\Div_{\partial B_R}v\:dS.
\end{eqnarray*}
We insert this expression into \eqref{intermdiate3}. We find
\begin{eqnarray*}
\oint_{\partial B_R}\ell_2\:dS
&=&
-
\oint_{\partial B_R}(\delta_{i}v_{i})^2\:dS
-
\frac{n-2}{R}\oint_{\partial B_R}v_j\:\delta_{j}v_i\nu_{i}\:dS\\
\nonumber&&+
\frac{n-2}{R}\oint_{\partial B_R}N\:\Div_{\partial B_R}v\:dS.
\end{eqnarray*}
Thus,
\begin{eqnarray*}
\oint_{\partial B_R}(\ell_2 +
(\p_i^*\tilde v\cdot x_{\xi_i})^2)
\:dS&=&
-
\frac{n-2}{R}\oint_{\partial B_R}v_j\:\delta_{j}v_i\nu_{i}\:dS\\
&&+
\frac{n-2}{R}\oint_{\partial B_R}N\:\Div_{\partial B_R}v\:dS.
\end{eqnarray*}
This identity, together with \eqref{intermed1}, and the fact that
\begin{eqnarray*}
N\:\Div_{B_R}v -v_{j}\:\delta_{j}v_{i}\:\nu_{i}=N\:\Div v -v_{j}\:\partial_{j}v_{i}\:\nu_{i}
\end{eqnarray*}
implies the following lemma.
%%%%%%%%%%%%%%%%%%%%%%%%%%%%%%%%%
\begin{lemma}\label{mball2}
For an arbitrary vector field $v=v^\tau +N\nu$ the second variation $\ddot{S}_{B_R}(0)$ assumes the form
\begin{eqnarray*}\oint_{\p B_R} \ddot{m}(0)\:dS& =&\oint_{\partial B_R}\left(|\nabla^\tau N|^2 -\frac{(n-1)}{R^2} N^2\right)\:dS\\
&&+
\frac{n-1}{R}\oint_{\partial B_R}\left(N\:\Div v -v_{j}\:\partial_{j}v_{i}\:\nu_{i}+w\cdot\nu\right)\:dS.
\end{eqnarray*}
\end{lemma}
%%%%%%%%%%%%%%%%%%%%%%%%%%%%%%%%%
Let us now consider vector fields, which are volume preserving of the second order (Lemma \ref{volpres1}). We observe that, in view of \eqref{volcomp1}, the second integral on the right-hand side in Lemma \ref{mball2} vanishes. Therefore 
\begin{eqnarray}\label{isodefrad1}
&&\oint_{\partial  B_R}\ddot{m}(0)\:dS
=
\oint_{\partial B_R}|\nabla^\tau N|^2\:dS
-
\frac{n-1}{R^2}\oint_{\partial B_R}N^2 \:dS.
\end{eqnarray}
%%%%%%%%%%%%%%%%%%%%%%%%%%%%%RRRRRRRRRRRRRRR
\begin{remark}
It is interesting to observe that the second variation $\oint_{\partial  B_R}\ddot{m}(0)\:dS$ neither depends on $w$ nor on the tangential components of the vector field $v$. Formula \eqref{isodefrad1} is very similar to some results in \cite{Li12} in connection with stability inequalities for immersed submanifolds.
\end{remark}
%%%%%%%%%%%%%%%%%%%%rrrrrrrrrrrrrrrrrrrrrrrrrr
Let us introduce the following notation.
\begin{eqnarray*}
{\cal{S}}(t) := \oint_{\p B_R} m(t)\:dS \qquad(=\hbox{surface area of $\p \Omega_t$}).
\end{eqnarray*}
Next we determine all volume preserving vector fields of first and second order for which the second variation
$\ddot{\cal{S}}(0)$  vanishes. They form {\sl the kernel} of $\ddot{\cal{S}}(0)$. For this purpose we recall the eigenvalue problem 
\begin{align}\label{Stek}
\Delta_{\mathbb{S}^{n-1}} \phi + \mu \phi=0 \tx{on} \mathbb{S}^{n-1}.
\end{align}
It is well-known that the eigenfunctions are the spherical harmonics of order $k$. The corresponding eigenvalues are $\mu_k=k(k+n-2)$, $k\in \mathbb{N}^+$ with the multiplicity
$(2k+n-2)\frac{(k+n-3)!}{(n-2)!k!}$. If $v$ is volume preserving of the first order, then $ \oint_{\p B_R}N\:dS=0$. By the variational characterization of the eigenvalues we get 
\begin{align*}
\oint_{\partial B_R}|\nabla^\tau N|^2\:dS
\geq \frac{\mu_2}{R^2}\oint_{\partial B_R}N^2 \:dS=\frac{n-1}{R^2}\oint_{\partial B_R}N^2 \:dS.
\end{align*}
Equality holds if and only if $(v\cdot \nu)$ is an element of the eigenspace corresponding to $\mu= (n-1)/R^2$
or if $N=0$.
%%%%%%%%%%%%%%%%%%%%
\begin{example}
Suppose that on $\p B_R$ the vector field $v$ points only in tangential direction. Then
$v= g^{ij}(v\cdot x_{\xi_j})x_{\xi_i}$. The vector $y=x(\xi')+t\tilde v(\xi')$ is orthogonal to $\p B_R$. Its length is
$|y|^2= R^2+t^2g^{ij}(v\cdot x_{\xi_i})(v\cdot x_{\xi_j})$. The boundary $\p \Omega_t$ can therefore be represented by $y+ \frac{t^2}{2}g_0 \nu+o(t^2)$, where $g_0=g^{ij}(v\cdot x_{\xi_i})(v\cdot x_{\xi_j})$ and $\nu =\frac{y}{R}$. The domain $\Omega_t$ is therefore a second order perturbation of $B_R$. 
\end{example}
%%%%%%%%%%%%%%%%%%%%LLLLLLLLLLL
\begin{definition} 
\sl{A perturbation of the form
\begin{align*}
y=x+ tN\nu +\frac{t^2}{2}(w\cdot\nu)\nu +o(t^2) 
\end{align*}
on $\partial B_R$ is called a Hadamard perturbation. }
\end{definition}
%%%%%%%%%%%%%%%%%%%%LLLLLLLLLLL
From the previous consideration it follows immediately that every small perturbation of the ball can be described by a Hadamard perturbation. Consequently we have
%%%%%%%%%%%%%%%%%%%%LLLLLLLLLLL
\begin{lemma}\label{Hadamard} Assume $N \neq a_ix_i$ on $\p B_R$. Then $\ddot{\cal{S}}(0)>0$ for every Hadamard perturbation.
\end{lemma}
%%%%%%%%%%%%%%%%%%%%%%%%%%%%%%%%%%%%%%
%%%%%%%%%%%%%%%%%%%%%%%%%%%%%%%%%%%%%
\section{Energies}\label{energ}
%%%%%%%%%%%%%%%%%%%%%%%%%%%%%%%%%%%%%%
%%%%%%%%%%%%%%%%%%%%%%%%%%%%%%%%%%%%%%
Let $(\Omega_t)_t$ be a family of domains described in the previous chapter, and let  $G:\rz\to\rz$ denote a smooth function, i.e. $G\in C^2_{loc}(\rz)$ at least. We denote by $g$ its derivative: $G'=g$. Consider the energy functional
\begin{eqnarray}\label{rt}
{\cal{E}}(\Omega_t,u)&:=& \int_{\Omega_t}\vert \nabla_{y}u\vert^2\:dy-2\int_{\Omega_t}G(u)\:dy+\alpha\oint_{\p\Omega_t}u^2\:dS_t.
\end{eqnarray}
The boundary integral represents the elastic energy. We shall call ${\cal{E}}(\Omega_t,u)$ the {\sl Robin energy}.
\newline
\newline
A critical point $\tilde{u}\in H^{1}(\Omega_t)$ of \eqref{rt} satisfies the Euler Lagrange equation
\begin{eqnarray}
\label{rt01}\Delta_{y}\tilde{u}+g(\tilde{u})&=&0\qquad\hbox{in}\:\:\Omega_t\\
\label{rt02}\partial_{\nu_t}\tilde{u}+\alpha\tilde{u}&=&0\qquad\hbox{in}\:\:\partial\Omega_t,
\end{eqnarray}
where $\nu_t$  stands for the outer normal of $\partial\Omega_t$. A special case is the torsion problem with $G(w)=w$.
\newline
\newline
Assume that $\tilde{u}$ solves \eqref{rt01} - \eqref{rt02}. We set
\begin{eqnarray*}
{\cal{E}}(t):={\cal{E}}(\Omega_t,\tilde{u}).
\end{eqnarray*}
In a first step we transform the integrals onto $\Omega$ and $\partial\Omega$. Let 
\begin{eqnarray*}
y=x+tv(x)+\frac{t^2}{2}w(x)+o(t^2)
\end{eqnarray*} 
be defined as in \eqref{yvar} and let $x(y)$ be its inverse.
After a change of variables we get
\begin{eqnarray}\label{transen}
{\cal{E}}(t)&=& \int_{\Omega}\partial_{i}\tilde{u}(t)\:\partial_{j}\tilde{u}(t)\left(\frac{\partial x_i}{\partial y_k}\right)\left(\frac{\partial x_j}{\partial y_k}\right)\:J(t)\:dx\\
&&\nonumber-2\int_{\Omega}G(\tilde{u}(t))\:J(t)\:dx\\
&&\nonumber
+\alpha\oint_{\p\Omega}\tilde{u}^2(t)\:m(t)\:dS,
\end{eqnarray}
where $\tilde{u}(t):=\tilde{u}(x+tv(x)+\frac{t^2}{2}w(x),t),\:t\in (-\epsilon,\epsilon)$.
We set
\begin{eqnarray}\label{aij1}
A_{ij}(t):=\frac{\partial x_i}{\partial y_k}\frac{\partial x_j}{\partial y_k}\:J(t).
\end{eqnarray}
The expression \eqref{transen} assumes now the concise form
\begin{eqnarray}\label{transeng}
{\cal{E}}(t)= \int_\Omega \nabla \tilde u A \nabla\tilde  u\:dx -2\int_\Omega G(\tilde u)J\:dx+\alpha \oint_{\p \Omega} \tilde u^2m\:dx.
\end{eqnarray}
Thus, in the domain $\Omega$ the solution $\tilde{u}(t)$ solves the transformed equation
\begin{eqnarray}
\label{transsol1g}L_A \tilde{u}(t)+g(\tilde{u}(t))\:J(t)&=&0\qquad\hbox{in}\:\:\Omega\\
\label{transsol2g}\partial_{\nu_A}\tilde{u}(t)+\alpha\:m(t)\:\tilde{u}(t)&=&0\qquad\hbox{in}\:\:\partial\Omega,
\end{eqnarray}
where
\begin{eqnarray}\label{transcoeff}
L_A=\partial_{j}\left(A_{ij}(t)\partial_{i}\right)\qquad\qquad\hbox{and}\qquad
\qquad\partial_{\nu_A} =\nu_{i}A_{ij}(t)\partial_{j}.
\end{eqnarray}
It is convenient to write the equations \eqref{transsol1g} - \eqref{transsol2g} for $\tilde u$ in the weak form
\begin{align}\label{Euler1}
\int_\Omega \nabla \phi A\nabla \tilde u\:dx +\alpha \oint_{\p \Omega} \phi\: \tilde u \:m\:dS= \int_\Omega g(\tilde u)\:\phi \:J\:dx, \: \forall \phi \in W^{1,2}(\Omega).
\end{align}
\newline
\newline
We now consider the eigenvalue problem
\begin{eqnarray}
\label{re1}\Delta_{y}\tilde{u}+\lambda(\Omega_t)\tilde{u}&=&0\qquad\hbox{in}\:\:\Omega_t\\
\label{re2}\partial_{\nu_t}\tilde{u}+\alpha\tilde{u}&=&0\qquad\hbox{in}\:\:\partial\Omega_t.
\end{eqnarray}
The eigenvalue $\lambda(t):=\lambda(\Omega_t)$ in \eqref{re1} and \eqref{re2} is characterized by the Rayleigh quotient
\begin{eqnarray}\label{re}
\lambda(t)={\cal{R}}(t)&:=&\frac{\int_{\Omega_t}\vert \nabla_{y}\tilde{u}\vert^2\:dy+\alpha\oint_{\p\Omega_t}\tilde{u}^2\:dS}{\int_{\Omega_t}\tilde{u}^2\:dy}.
\end{eqnarray}
The change of variable \eqref{yvar} yields
\begin{eqnarray}\label{transew}
{\cal{R}}(t)=\frac{\int_{\Omega}A_{ij}(t)\:\partial_{i}\tilde{u}(t)\:\partial_{j}\tilde{u}(t)\:dx+\alpha\oint_{\p\Omega}\tilde{u}(t)^2\:m(t)\:dS}{\int_{\Omega}\tilde{u}^2(t)\:J(t)\:dx}.
\end{eqnarray}
Thus, in the domain $\Omega$ the solution $\tilde{u}(t)$ solves the transformed equation
\begin{eqnarray}
\label{transew1}L_A \tilde{u}(t)+\lambda(t)J(t)\:\tilde{u}(t)&=&0\qquad\hbox{in}\:\:\Omega\\
\label{transew2}\partial_{\nu_A}\tilde{u}(t)+\alpha\:m(t)\:\tilde{u}(t)&=&0\qquad\hbox{in}\:\:\partial\Omega.
\end{eqnarray}
We use $\tilde u$ as a test function in the above equation We obtain the identity
\begin{eqnarray}\label{eigen123}
\int_\Omega \nabla \tilde u A \nabla \tilde u\:dx +\alpha \oint_{\p \Omega} \tilde u^2m \: dS =\lambda (t)\int_\Omega \tilde u^2 J\:dx,
\end{eqnarray}
which will be used later.
%%%%%%%%%%%%%%%%%%%%%%%%%%%%%%%%%%%%%%%%%
\subsection{Expansions}
\label{enexp}
%%%%%%%%%%%%%%%%%%%%%%%%%%%%%%%%%%%%%%%%%
In this subsection we assume that \eqref{rt01} - \eqref{rt02} has a unique solution. 
We expand formally all relevant quantities with respect to $t$. Under suitable regularity assumption on $\Omega_t$ such processes can be justified.
\newline
\newline
We start with the energy \eqref{transeng}
\begin{eqnarray*}
{\cal{E}}(t)= \int_\Omega \nabla \tilde u A \nabla\tilde  u\:dx -2\int_\Omega G(\tilde u)J\:dx+\alpha \oint_{\p \Omega} \tilde u^2m\:dx,
\end{eqnarray*}
where $\tilde u$ is a weak solution of \eqref{Euler1}. 
\newline
\newline
Recall that $\tilde{u}(t)=\tilde{u}(x+tv(x)+\frac{t^2}{2}w(x),t),\: t\in ]-\epsilon,\epsilon[$. Under sufficient regularity the following expansion is valid
\begin{eqnarray*}
\tilde{u}(t)=\tilde{u}(0)+t\:\dot{\tilde{u}}(0)+\frac{t^2}{2}\:\ddot{\tilde{u}}(0)+o(t^2).
\end{eqnarray*}
We set $u'(x):=\partial_{t}\tilde{u}(x+tv(x)+\frac{t^2}{2}w(x),t)\vert_{t=0}$ and get the following formulas for the coefficients of this expansion:
\begin{eqnarray}
\label{utilde1}\tilde{u}(0)&=&u(x)\qquad\hbox{and}\\
\label{utilde2}\nonumber\dot{\tilde{u}}(0)&=&\partial_{t}\tilde{u}(x+tv(x)+\frac{t^2}{2}w(x),t)\vert_{t=0}\\
\nonumber&+&(v(x)+t w(x))\cdot\nabla\tilde{u}(x+tv(x)+\frac{t^2}{2}w(x),t)\vert_{t=0}\\
&=&u'(x)+v(x)\cdot\nabla u(x).
\end{eqnarray}
The functions $\dot{\tilde{u}}(0)$ and $u'(x)$ are often called {\sl material} and {\sl shape derivative}, respectively.
\newline
\newline
We also expand $A_{ij}(t)$ with respect to $t$:
\begin{eqnarray}\label{aij2}
A_{ij}(t)=A_{ij}(0)+t\:\dot{A}_{ij}(0)+\frac{t^2}{2}\ddot{A}_{ij}(0)+o(t^2).
\end{eqnarray}
Later we will compute $\dot{\cal{E}}(0)$ and $\ddot{\cal{E}}(0)$. For this purpose
 we shall need the explicit terms in \eqref{aij2}.
A lengthy but straightforward computation gives
%%%%%%%%%%%%%%%%%%%%%%%%%%%%%%%%%%%%%%%
\begin{lemma}\label{aij3}
\begin{eqnarray*}
A_{ij}(0)&=&\delta_{ij};\\
\dot{A}_{ij}(0)&=&\Div v\:\delta_{ij}-\partial_{j}v_{i}-\partial_{i}v_{j};\\
\ddot{A}_{ij}(0)&=&\left((\Div v)^2-D_v:D_v\right)\:\delta_{ij} +2\left(\partial_{k}v_i\:\partial_{j}v_k+\partial_{k}v_j\:\partial_{i}v_k\right)\\
&&
+
2\:\partial_{k}v_i\:\partial_{k}v_j
-
2\:\Div v\left(\partial_{j}v_{i}+\partial_{i}v_{j}\right)
+\Div\:w\: \delta_{ij}\\
&&
 -\p_iw_j-\p_j w_i.
\end{eqnarray*}
\end{lemma}
%%%%%%%%%%%%%%%%%%%%%%%%%%%%%%%%%%%%%%%%%
Finally we recall from \eqref{jform}
\begin{eqnarray}
\label{jform1}J(0)&=& 1\\
\label{jform2}\dot{J}(0)&=&\Div v\\
\label{jform3}\ddot{J}(0)&=&\left(\Div v\right)^2-D_v:D_v+\Div\:w.
\end{eqnarray}
%%%%%%%%%%%%%%%%%%%%%%%%%%%%%%%%%%%%%%%%%
\subsection{Differentiation of the energy and the eigenvalue}
\label{diffenei}
%%%%%%%%%%%%%%%%%%%%%%%%%%%%%%%%%%%%%%%%%
\subsubsection{First and Second Variation}
The first derivative of ${\cal{E}}(t)$ with respect to the parameter $t$ is called {\sl{first domain variation}}, and the second derivative is called the {\sl{second domain variation}}. In modern texts often the expression {\sl{shape gradient}} and {\sl{shape Hessian}} are used instead.
\label{subsubsv}
Direct computation gives
\begin{eqnarray}\label{eliminate}
\dot{\cal{E}}(t)&=&\int_\Omega \nabla \tilde u \dot A \nabla \tilde u \:dx-2\int_\Omega G(\tilde u) \dot J \:dx +\alpha \oint_{\p \Omega} \tilde u^2 \dot m\:dS\\
\nonumber&&+2\int_\Omega \nabla \dot{\tilde u} A \nabla \tilde u\:dx
-2\int_\Omega g\dot{\tilde u}J\:dx +2\alpha \oint_{\p \Omega} \dot{\tilde u}\tilde u m\:dS.
\end{eqnarray}
In \eqref{Euler1} we choose $\phi=\dot{\tilde u}$. We then can eliminate the terms containing $\dot{\tilde u}$ in \eqref{eliminate}. We obtain 
\begin{align}\label{edot}
\dot{\cal{E}}(t)=\int_\Omega\nabla \tilde u \dot A \nabla \tilde u \:dx + \alpha \oint_{\p \Omega} \tilde u^2 \dot m\:dS -2\int_\Omega G\dot J\:dx.
\end{align}
Note that $\dot{\cal{E}}(t)$ is independent of $\dot{\tilde u}$.
\newline
\newline
Next we want to find an expression for the second derivative.
Differentiation of \eqref{Euler1} implies 
\begin{eqnarray}\label{Euler2}
&&\int_\Omega\left[ \nabla \phi A \nabla \dot {\tilde u} + 
\nabla \phi \dot A \nabla \tilde u\right] \:dx\\
\nonumber&&\qquad +\alpha \oint_{\p \Omega}
( \phi \dot{ \tilde u} m + \phi\tilde u\dot m)\:dS
=\int_\Omega (g'(\tilde{ u})\dot {\tilde u} \phi J+
g(\tilde u)\phi\dot J)\:dS
\end{eqnarray}
for all $\phi \in W^{1,2}(\Omega)$.
\newline
Differentiation  of \eqref{edot} yields
\begin{eqnarray}\label{eddot}
\ddot{\cal{E}}(t)&=&\int_\Omega\left[ \nabla \tilde u \ddot A \nabla \tilde u +2 \nabla \dot{\tilde u} \dot A \nabla \tilde u- 2g\dot{\tilde u}\dot J -2G\ddot J\right]\:dx\\
\nonumber&&+\alpha \oint_{\p \Omega}(2\tilde u\dot{\tilde u} \dot m + \tilde u^2\ddot m)\:dS.
\end{eqnarray}
By means of \eqref{Euler2} with $\phi= \dot{\tilde u}$ we get
\begin{eqnarray}\label{eddot2}
\ddot{\cal{E}}(t)&=&\int_\Omega \nabla \tilde u \ddot A \nabla \tilde u \:dx +\alpha \oint_{\p \Omega} \tilde u^2\ddot m\:dS \\
\nonumber&&-2\int_\Omega G\ddot J\:dx
-2\int_\Omega \nabla \dot{\tilde u} A \nabla \dot{\tilde u} \:dx -2\alpha \oint_{\p \Omega}\dot{\tilde u}^2m\:dS\\
\nonumber&&+2\int_\Omega g' \dot {\tilde u}^2J\:dx. 
\end{eqnarray}
In accordance with the first derivative which does not depend on $\dot{\tilde u}$, the second derivative does not depend on $\ddot{\tilde u}$.
\newline
\newline
In order to compute the variations of the eigenvalue 
we first recall that ${\tilde{u}}$ solves \eqref{eigen123}. We impose the normalization
\begin{eqnarray}\label{norma}
\int_{\Omega}\tilde{u}^2(t)\:J(t)\:dx =1.
\end{eqnarray}
This implies\begin{eqnarray}\label{eigennorm}
\frac{d}{dt}\int_{\Omega}\tilde{u}^2(t)\:J(t)\:dx =2\int_{\Omega}\tilde{u}(t)\:\dot{\tilde{u}}(t)\:J(t)\:dx+\int_{\Omega}\tilde{u}(t)^2\:\dot{J}(t)\:dx=0.
\end{eqnarray}
We differentiate \eqref{re} and use the normalization \eqref{norma}. Then we choose $\phi=\dot{\tilde{u}}$ as a test function in the weak formulation of \eqref{transew1} and \eqref{transew2}. By \eqref{eigennorm} we get
\begin{eqnarray}\label{dotlambda1}
\dot{\lambda}(t)
=
\int_\Omega \nabla \tilde u \dot A \nabla \tilde u\:dx 
+
\alpha \oint_{\p \Omega} \tilde u^2\dot m \: dS
-
\lambda(t) \int_\Omega \tilde u^2 \dot{J}\:dx.
\end{eqnarray}
Thus, $\dot{\lambda}(t)$ does not depend on $\dot{\tilde{u}}(t)$. 
\newline
\newline
We differentiate \eqref{dotlambda1} with respect to $t$. Then we differentiate \eqref{eigen123} with respect to $t$ and choose $\phi=2\dot{\tilde{u}}$ in the weak formulation of \eqref{transew1} and \eqref{transew2}. By \eqref{eigennorm}, and under the assumption $\dot{\lambda}(0)=0$, we get for $t=0$
\begin{eqnarray}\label{dotlambda2}
\ddot{\lambda}(0)
&=&
\int_{\Omega} \nabla \tilde u \ddot A \nabla \tilde u\:dx 
-
2\int_{\Omega}\vert \nabla \dot{\tilde{u}} \vert^2\:dx 
+
\alpha \oint_{\p \Omega} \tilde{u}^2\ddot m \: dS\\
\nonumber&&
-
2\alpha \oint_{\p \Omega} \dot{\tilde u}^2\:m \: dS-
\lambda(0) \int_\Omega \tilde u^2 \ddot{J}\:dx
+
2\lambda(0)\int_\Omega \dot{\tilde u}^2\:J\:dx .
\end{eqnarray}
Thus, $\ddot{\lambda}(0)$ does not depend on $\ddot{\tilde{u}}(0)$. 
%%%%%%%%%%%%%%%%%%%%%%%%%%%%%%%%%%%%%
\subsubsection{Third variation}
\label{subsubtv}
%%%%%%%%%%%%%%%%%%%%%%%%%%%%%%%%%%%%%
In order to compute the third variation $\dddot{\cal{E}}(s)$ we proceed exactly in the same way as before.
We differentiate \eqref{eddot2}.
\begin{eqnarray*}
\dddot{\cal{E}}(t)&=&
2\int_{\Omega}\nabla\dot{\tilde{u}}\cdot\left(\ddot{A}\nabla\tilde{u}\right)\:dx
+
2\int_{\Omega}\nabla\tilde{u}\cdot\left(\dddot{A}\nabla\tilde{u}\right)\:dx
+
2\alpha\oint_{\partial\Omega}\tilde{u}\:\dot{\tilde{u}}\:\ddot{m}\:dS\\
&&
+
\alpha\oint_{\partial\Omega}\tilde{u}^2\:\dddot{m}\:dS
-
2\int_{\Omega}G'(\tilde{u})\:\dot{\tilde{u}}\:\ddot{J}\:dx
-
2\int_{\Omega}G(\tilde{u})\:\dddot{J}\:dx\\
&&
-
4\int_{\Omega}\nabla\ddot{\tilde{u}}\cdot\left(A\nabla\dot{\tilde{u}}\right)\:dx
-
2\int_{\Omega}\nabla\dot{\tilde{u}}\cdot\left(\dot{A}\nabla\dot{\tilde{u}}\right)\:dx
-
4\alpha\oint_{\partial\Omega}\dot{\tilde{u}}\:\ddot{\tilde{u}}\:m\:dS\\
&&
-
2\alpha\oint_{\partial\Omega}\dot{\tilde{u}}^2\:\dot{m}\:dS
+
2\int_{\Omega}G'''(\tilde{u})\:\dot{\tilde{u}}^3\:J\:dx
+
4\int_{\Omega}G''(\tilde{u})\:\dot{\tilde{u}}\:\ddot{\tilde{u}}\:J\:dx\\
&&
+
2\int_{\Omega}G''(\tilde{u})\:\dot{\tilde{u}}^2\:\dot{J}\:dx.
\end{eqnarray*}
Differentiation of \eqref{Euler2} gives
\begin{eqnarray*}
&&\int_\Omega(
g'(\tilde u)\ddot {\tilde u}J\phi +g''(\tilde u)\dot {\tilde u}^2J\phi+2g'(\tilde u)\dot {\tilde u}\dot J\phi +g(\tilde u)\ddot J\phi\:dx\\
&=& \int_\Omega( \nabla \phi \ddot A \nabla \tilde u+2\nabla \phi \dot A \nabla \dot{\tilde u} +\nabla \phi A\nabla \ddot{\tilde u})\:dx \\
&&+ \alpha \oint_{\p\Omega}( \tilde u\phi  \ddot m +2\dot{\tilde u}\phi \dot m +\ddot{\tilde u} \phi m) \:dS 
\end{eqnarray*}
If $\phi=4\tilde{u}(t)$, then
\begin{eqnarray}
&&\nonumber-4\int_{\Omega}\nabla\ddot{\tilde{u}}\cdot\left(A\nabla\dot{\tilde{u}}\right)\:dx
+
4\oint_{\partial\Omega}\dot{\tilde{u}}\:\partial_{\nu_{A}}\ddot{\tilde{u}}\:m\:dS
+
4\int_{\Omega}G''(\tilde{u})\:\dot{\tilde{u}}\:\ddot{\tilde{u}}\:J\:dx\\
&&\label{weak3dot}\quad
-
4\int_{\Omega}\nabla\dot{\tilde{u}}\cdot\left(\ddot{A}\nabla\tilde{u}\right)\:dx
+
4\oint_{\partial\Omega}\dot{\tilde{u}}\:\partial_{\nu_{\ddot{A}}}\tilde{u}\:m\:dS\\
&&\nonumber\quad
-
8\int_{\Omega}\nabla\dot{\tilde{u}}\cdot\left(\dot{A}\nabla\dot{\tilde{u}}\right)\:dx
+
8\oint_{\partial\Omega}\dot{\tilde{u}}\:\partial_{\nu_{\dot{A}}}\dot{\tilde{u}}\:m\:dS
+
4\int_{\Omega}G'''(\tilde{u})\:\dot{\tilde{u}}^3\:J\:dx\\
&&\nonumber\quad+
8\int_{\Omega}G''(\tilde{u})\:\dot{\tilde{u}}^2\:\dot{J}\:dx
+
4\int_{\Omega}G'(\tilde{u})\:\dot{\tilde{u}}^2\:\ddot{J}\:dx=0.
\end{eqnarray}
Note that only three integrals in \eqref{weak3dot} contain $\ddot{\tilde{u}}$. They also appear in $\dddot{{\cal{E}}}(t)$. Thus, $\dddot{{\cal{E}}}(t)$ does not depend on  $\ddot{\tilde{u}}$. Hence
\begin{eqnarray}\label{thirdderen}
\nonumber\dddot{\cal{E}}(t)&=&
2\int_{\Omega}\nabla\dot{\tilde{u}}\cdot\left(\ddot{A}\nabla\tilde{u}\right)\:dx
+
2\int_{\Omega}\nabla\tilde{u}\cdot\left(\dddot{A}\nabla\tilde{u}\right)\:dx
+
8\int_{\Omega}\nabla\dot{\tilde{u}}\cdot\left(\dot{A}\nabla\dot{\tilde{u}}\right)\:dx\\
&&-
6\int_{\Omega}G'(\tilde{u})\:\dot{\tilde{u}}\:\ddot{J}\:dx
-
2\int_{\Omega}G(\tilde{u})\:\dddot{J}\:dx
-
6\int_{\Omega}G''(\tilde{u})\:\dot{\tilde{u}}^2\:\dot{J}\:dx\\
&&\nonumber
-
2\int_{\Omega}G'''(\tilde{u})\:\dot{\tilde{u}}^3\:J\:dx
+
6\alpha\oint_{\partial\Omega}\tilde{u}\:\dot{\tilde{u}}\:\ddot{m}\:dS
+
6\alpha\oint_{\partial\Omega}\dot{\tilde{u}}^2\:\dot{m}\:dS\\
&&\nonumber
+
\alpha\oint_{\partial\Omega}\tilde{u}^2\:\dddot{m}\:dS.
\end{eqnarray}
Similarly we compute the third variation of $\lambda$.
We differentiate $\ddot{\lambda}(t)$ with respect to $t$. Then we differentiate 
the weak formulation of \eqref{transew1} and \eqref{transew2} twice with respect to $t$ and choose $\phi=-4\dot{\tilde{u}}$. With \eqref{eigennorm} we get
\begin{eqnarray}\label{dotlambda3}
\dddot{\lambda}(t)
&=&
\int_{\Omega} \nabla \tilde u \dddot A \nabla \tilde u\:dx 
+
6\int_{\Omega} \nabla \dot{\tilde{u}} \ddot{A} \nabla\tilde{u}\:dx 
+
6\int_{\Omega} \nabla \dot{\tilde{u}} \dot{A} \nabla\dot{\tilde{u}}\:dx\\
\nonumber&&+
\alpha \oint_{\p \Omega} \tilde{u}^2\dddot m \: dS
+
6\alpha \oint_{\p \Omega} \tilde u\:\dot{\tilde u}\:\ddot{m} \: dS
+
6\alpha \oint_{\p \Omega} \dot{\tilde u}^2\:\dot{m} \: dS\\
\nonumber&&-
\lambda(t) \int_\Omega \tilde u^2 \dddot{J}\:dx
-
6\lambda(t)\int_\Omega \tilde u\:\dot{\tilde u}^2\:\ddot{J}\:dx
-
3\dot{\lambda}(t)\int_\Omega {\tilde u}^2\:\ddot{J}\:dx\\
\nonumber&&
-
6\dot{\lambda}(t)\int_\Omega \dot{\tilde u}^2\:J\:dx
-
6\lambda(t)\int_\Omega \dot{\tilde u}^2\:\dot{J}\:dx
-
12\dot{\lambda}(t)\int_\Omega \tilde u\:\dot{\tilde u}\:\dot{J}\:dx.
\end{eqnarray}
Thus, $\dddot{\lambda}(t)$ does not depend on $\ddot{\tilde{u}}$.
A direct consequence is
%%%%%%%%%%%%%%%%%%%%%%%%%%%CCCCCCCCCCCC
\begin{corollary} The derivatives of ${\cal{E}}(t)$ and of $\lambda(t)$ of order greater than two are expressed in terms of the derivatives of $\tilde u$ of  two orders lower.
\end{corollary}
%%%%%%%%%%%%%%%%%%%%%%%%%%%%%%cccccccccccccc
This phenomenon was observed by D. D. Joseph \cite{Jo67} for the eigenvalues.
%%%%%%%%%%%%%%%%%%%%%%%%%%%%%%%%%%%%%%
%%%%%%%%%%%%%%%%%%%%%%%%%%%%%%%%%%%%%
\section{The first variation}
\label{firstvar}
%%%%%%%%%%%%%%%%%%%%%%%%%%%%%%%%%%%%%
%%%%%%%%%%%%%%%%%%%%%%%%%%%%%%%%%%%%%%
\subsection{Energies}
\label{energy}
The goal of this section is to represent $\dot{\cal{E}}(0)$ as a boundary integral. By \eqref{edot} we have
\begin{eqnarray*}
\dot{\cal{E}}(0)=\underbrace{ \int_\Omega \partial_{i} u\dot A_{ij}(0)\partial_{j} u \:dx}_{\dot {\cal{E}}_1} +\alpha \oint_{\p \Omega}  u^2 \dot m(0)\:dS -2\underbrace{\int_\Omega G(u)\dot J(0)\:dx}_{\dot {\cal{E}}_2}.
\end{eqnarray*}
From Lemma \ref{aij3} we conclude, after integration by parts, that
\begin{equation}
\dot {\cal{E}}_1= \oint_{\p \Omega}\{|\nabla u|^2(v\cdot \nu)-2(v\cdot \nabla u)(\nu\cdot \nabla u)\}\:dS -2\int_\Omega g(u)(v\cdot \nabla u)\:dx.
\end{equation}
Moreover,
\begin{equation}
\dot {\cal{E}}_2= \oint_{\p \Omega} G(u)(v\cdot \nu)\:dS -\int_\Omega g(u)(v\cdot \nabla u)\:dx.
\end{equation}
Hence by \eqref{SpurA} and the boundary condition  \eqref{rt02} for $u$
\begin{eqnarray*}
\dot{\cal{E}}(0)&=&\oint_{\p \Omega}\{|\nabla u|^2-2G(u))(v\cdot \nu) +2\alpha (v\cdot \nabla u)u \\
&&\qquad+\alpha u^2(\Div_{\p \Omega}v^\tau +(n-1)(v\cdot \nu) H)\}\:dS.
\end{eqnarray*}
Observe that
\begin{equation}
v\cdot \nabla u= (v^\tau +(v\cdot \nu)\nu)\cdot (\nabla^\tau u +(\nu\cdot \nabla u)\nu)= v^\tau\cdot \nabla^\tau u -\alpha (v\cdot \nu)u.
\end{equation}
Thus,
\begin{eqnarray*}
\dot{\cal{E}}(0)&=&\oint_{\p \Omega}(v\cdot \nu) 
\{|\nabla u|^2-2G(u) -2\alpha^2 u^2 +\alpha (n-1)Hu^2\}\:dS\\
&&\qquad+\alpha \oint_{\p \Omega}(2v^\tau u\nabla^\tau u +u^2 \Div_{\p \Omega} v^\tau)\:dS.
\end{eqnarray*}
The last integral vanishes by \eqref{tanpi}. Finally we have 
\begin{align}\label{Edot1}
\dot{\cal{E}}(0)=\oint_{\p \Omega}(v\cdot \nu) 
\{|\nabla u|^2-2G(u) -2\alpha^2 u^2 +\alpha (n-1)Hu^2\}\:dS.
\end{align}
In particular we observe that $\dot{\cal{E}}(0)=0$ for all purely tangential deformations.
From the expression \eqref{Edot1} above we deduce
\newline
%%%%%%%%%%%%%%%%%%%%%%%%%%%%%%%%%%%%%%%%%
\begin{theorem}\label{overdet1}
Let $\Omega_t$ be a family of volume preserving perturbations  of $\Omega$ as described in \eqref{yvar}. Then $\Omega$ is a critical point of the energy ${\cal{E}}(t)$, i.e. $\dot{{\cal{E}}}(0)=0$
if and only if
\begin{eqnarray}\label{overrt}
\vert\nabla u\vert^2-2G(u)-2\alpha^2 u^2+\alpha(n-1)u^2H=const.\qquad\hbox{on}\quad\partial\Omega.
\end{eqnarray}
\end{theorem}
%%%%%%%%%%%%%%%%%%%%%%%%%%%%%%%%%%%%%%%%%
{\bf{Proof}} 
Write for short
\begin{equation}
z(x):=\vert\nabla u\vert^2-2G(u)-2\alpha^2 u^2+\alpha(n-1)u^2H \tx{and} \overline{z}:=|\p\Omega|^{-1}\oint_{\p \Omega} z\:dS.
\end{equation}
Then, since $\oint_{\p\Omega}(v\cdot \nu)\:dS=0$,
\begin{equation}
\oint_{\p\Omega} (v\cdot \nu) z\:dS= \oint_{\p\Omega} (v\cdot \nu) (z-\overline{z})\:dS.
\end{equation}
Put $Z^{\pm} =$max$\{0,\pm (z-\overline{z})\}$. Hence
\begin{equation}
\oint_{\p\Omega} (v\cdot \nu) z\:dS=\oint_{\p\Omega} (v\cdot \nu)(Z^+-Z^-)\:dS.
\end{equation}
Suppose that $z\neq$const. Then $Z^\pm \neq 0$ and we can construct a volume
preserving perturbation such that $(v\cdot \nu)>0$ in supp$Z^+$ and $(v\cdot \nu)<0$ in supp$Z^-$. In this case we get $\dot{{\cal{E}}}(0)>0$, which is obviously a contradiction.
\hfill{$\rlap{$\sqcap$}\sqcup$}
\newline\newline 
%%%%%%%%%%%%%%%%%%%%%%%%%%%%%%%%%%%%%%%%
\begin{example} If $\Omega =B_R$ and $u(x)=u(|x|)$, then $\dot{\cal{E}}(0)=0$.
The question arises: are there  domains other than the ball for which we can find 
a solution $u:\Omega\to\mathbb{R}$ of the overdetermined problem
\begin{eqnarray*}
\Delta u+g(u)&=&0\qquad\hbox{in}\:\Omega\\
\partial_{\nu}u+\alpha u&=&0\qquad\hbox{in}\:\partial\Omega\\
\vert\nabla u\vert^2-2G(u)-2\alpha ^2u^2+\alpha(n-1)u^2H&=&
const.\tx{in}\partial\Omega?
\end{eqnarray*}
The moving plane method proposed by Serrin in \cite{Se71} doesn't seem to apply here.
\end{example}
%%%%%%%%%%%%%%%%%%%%%%%%%%%%%%%%%%%%%%%%%%%
\subsection{Eigenvalues}
\label{subeig}
%%%%%%%%%%%%%%%%%%%%%%%%%%%%%%
The same arguments as in Section \ref{energy} imply that
\begin{align}\label{eigenvaluev1}
\dot \lambda (0)= \oint_{\p \Omega}(|\nabla u|^2-\lambda(0)u^2-2\alpha^2 u^2+\alpha (n-1)Hu^2)(v\cdot \nu)\:dS.
\end{align}
In analogy to Theorem \ref{overdet1} we get by the same arguments
\newline
%%%%%%%%%%%%%%%%%%%%%%%%%%%%%%%%%%%%%%%%%%%
\begin{theorem}\label{overdet2}
Let $\Omega_t$ be a family of volume preserving perturbations of $\Omega$ as described in \eqref{yvar}. Then $\Omega$ is a critical point of the principal eigenvalue $\lambda(t)$, i.e. $\dot{\lambda}(0)=0$
if and only if
\begin{eqnarray}\label{overrt2}
\vert\nabla u\vert^2-\lambda u^2-2\alpha^2 u^2+\alpha(n-1)u^2H
=const.\qquad\hbox{in}\quad\partial\Omega.
\end{eqnarray}
\end{theorem}
%%%%%%%%%%%%%%%%%%%%%%%%%%%%%%%%%%%%%%%%%%%
Let us now determine the sign of the constant in \eqref{overrt2} for the ball $B_R$.  We set $z=\frac{u_r}{u}$ and observe that
\begin{eqnarray*}
\frac{dz}{dr} +z^2 +\frac{n-1}{r}{z}+ \lambda=0  \tx{in} ]0,R[.
\end{eqnarray*}
At the endpoint
\begin{eqnarray*}
\frac{dz}{dr}(R) +\alpha^2 -\frac{n-1}{R} \alpha + \lambda =0.
\end{eqnarray*} 
We know that $z(0)=0$ and $z(R)= -\alpha$. Assume $\alpha >0$.  If $z_r(R) >0$, then there exists a number 
$\rho \in ]0,R[$ such that $z_r(\rho) =0$, $z(\rho)<0$ and $z_{rr}(\rho) \geq 0$. From the equation we get 
$z_{rr}(\rho) =\frac{n-1}{\rho^2}z <0$, which leads to a contradiction. Consequently
\begin{align}\label{signk}
A:=  -\alpha^2 +\frac{ n-1}{R}\alpha-\lambda\leq 0.
\end{align}
Similarly we prove that $A\geq 0$ if $\alpha <0$. 
Consequently we have for $\alpha >0 \:(<0)$
\begin{align}
\dot \lambda (0)<0 \:(>0)
\end{align}
for all volume increasing (decreasing) perturbations 
$\oint_{\p B_R}v\cdot \nu \:dS >0(<0)$. Note that this observation extends partly the result of Giorgi and Smits \cite{Gi} who proved that $\lambda(\Omega)>\lambda(B_R)$ for any $\Omega \subset B_R$. The result for negative $\alpha$ was observed in \cite{Ba77}.
%%%%%%%%%%%%%%%%%%%%%%%%%%%%%%%%%%%%%%
%%%%%%%%%%%%%%%%%%%%%%%%%%%%%%%%%%%%%
\section{An equation for $u'$}
\label{uprime}
%%%%%%%%%%%%%%%%%%%%%%%%%%%%%%%%%%%%%
%%%%%%%%%%%%%%%%%%%%%%%%%%%%%%%%%%%%%%
In this section we derive a boundary value problem for the function $u'$ defined in \eqref{utilde2}.
Let $\tilde{u}(t)$ solve \eqref{transsol1g} - \eqref{transsol2g}. If we differentiate with respect to $t$ and evaluate the derivative at $t=0$ we get
\begin{eqnarray}
\label{exprt1}L_{A(0)} \dot{\tilde{u}}(0)+L_ {\dot{A}(0)}\tilde{u}(0)\\
\nonumber+
g'(\tilde{u}(0))\dot{\tilde{u}}(0)\:J(0)
+
g(\tilde{u}(0))\dot{J}(0)&=&0\qquad\hbox{in}\:\:\Omega\\
\label{exprt2}\partial_{\nu_{A(0)}}\dot{\tilde{u}}(0)+ \partial_{\nu_{\dot{A}(0)}}\tilde{u}(0)+\alpha\:m(0)\:\dot{\tilde{u}}(0)\\
\nonumber+\alpha\:\dot{m}(0)\:\tilde{u}(0)&=&0\qquad\hbox{in}\:\:\partial\Omega.
\end{eqnarray}
From Lemma \ref{aij3} we then get $\Delta u'+g'(u)\:u'=0$ in $\Omega$.
\newline 
\newline
The computation for the boundary condition for $u'$ is more involved.
\begin{eqnarray*}
\partial_{\nu_{A(0)}}\dot{\tilde{u}}(0)
&=&\partial_{\nu}(v\cdot\nabla u)+\partial_{\nu}u'\\
\partial_{\nu_{\dot{A}(0)}}\tilde{u}(0)&=&\Div v\:\partial_{\nu}u-\nu_{j}\:\partial_{j}v_{i}\:\partial_{i}u-\partial_{i}u\:\partial_{i} v_{j}\:\nu_{j}\\
&=& \Div v\: \p_\nu u-\nu\cdot D_v\nabla u -\nabla u \cdot D_v\nu\\
\alpha\:m(0)\:\dot{\tilde{u}}(0)&=&\alpha\:v\cdot\nabla u+\alpha\:u'\\
\alpha\:\dot{m}(0)\:\tilde{u}(0)&=&\alpha(n-1)(v\cdot\nu)H\:u+\alpha\: u\:\Div_{\partial\Omega}v^\tau.
\end{eqnarray*}
We insert these expressions into \eqref{exprt2}. We then take into account \eqref{SpurA} and the boundary condition $\p  u_\nu + \alpha u=0$.  We obtain
\begin{eqnarray*}
\partial_{\nu}u'+\alpha u'&=& -\p_\nu(v\cdot\nabla u)+ \nabla  u\cdot D_v\nu+\nu \cdot D_v\nabla u-\alpha v\cdot\nabla u\\
&&+\alpha u(\Div v-(n-1)(v\cdot \nu)H -\Div_{\p \Omega}v^\tau).
\end{eqnarray*}
We observe that since $\Div_{\partial\Omega}\nu=(n-1)H$, 
\begin{eqnarray*}
\Div v=\Div_{\partial\Omega}v^\tau+(n-1)(v\cdot\nu)\:H+\nu_{i}\:\partial_{i}v_{j}
\:\nu_{j}\quad\hbox{on}\:\partial\Omega.
\end{eqnarray*}
Thus,
\begin{eqnarray*}
\partial_{\nu}u'+\alpha u'&=& -\p_\nu( v\cdot\nabla u)+ \nabla  u\cdot D_v\nu+
\nu \cdot D_v\nabla u\\
&&-\alpha \:v\cdot\nabla u+\alpha\: u\:\nu\cdot D_v\nu.
\end{eqnarray*}
In view of \eqref{tgradient} and the boundary condition for $u$ we have
\begin{eqnarray*}
\nabla u\cdot D_{v}\nu=-\alpha\:u\:\nu\cdot D_v\nu+\nabla^\tau u\cdot D_{v}\nu.
\end{eqnarray*}
Hence
\begin{eqnarray*}
\partial_{\nu}u'+\alpha u'= -\p_\nu( v\cdot\nabla u)+ \nabla^\tau  u\cdot D_v\nu
+\nu \cdot D_v\nabla u-\alpha \:v\cdot\nabla u.
\end{eqnarray*}
Thus, $u'$ solves
\begin{eqnarray}
\label{bcup1}\Delta u'+g'(u)\:u'&=&0\tx{in}\Omega\\
\label{bcup2}\partial_{\nu}u'+\alpha u'&=& -\p_\nu( v\cdot\nabla u)+ 
\nabla^\tau  u\cdot D_v\nu+\nu \cdot D_v\nabla u\\
\nonumber&&\qquad-\alpha \:v\cdot\nabla u\tx{in}\partial\Omega.
\end{eqnarray}
Analogously we get for the eigenvalue problem 
\begin{eqnarray}
\label{uprew}\Delta u'+\lambda(0) u' &+& \dot{\lambda}(0)u=0\tx{in}\Omega\\
\partial_{\nu}u'+\alpha u'&=& -\p_\nu( v\cdot\nabla u)+ \nabla^\tau  u\cdot D_v\nu
+\nu \cdot D_v\nabla u\\
\nonumber&&\qquad
-\alpha \:v\cdot\nabla u\tx{in}\partial\Omega.
\end{eqnarray}
%%%%%%%%%%%%%%%%%%%%%%%%%%%%%%%%%%%%%%%%%%
\begin{example}
1. Of special interest will be the case where $\Omega$ is the ball $B_R$ of radius $R$ centered at the origin, and $u$  is a radial solution of $\Delta u +g(u)=0$ in $B_R$ with $\p_\nu u + \alpha u=0$ on $\p B_R$. Then \eqref{bcup2} becomes
\begin{align}\label{ballrt0}
\partial_{\nu}u'+\alpha u'=\left (g(u(R))-\frac{\alpha(n-1)}{R}u(R)+\alpha^2 u(R)\right)(v\cdot \nu).
\end{align}
For the torsion problem $g(u)=1$ we have
\begin{eqnarray}\label{ballrt}
u(x)=\frac{R}{\alpha\:n}+\frac{1}{2n}\left(R^2-\vert x\vert^2\right).
\end{eqnarray}
We insert $u(R)=\frac{R}{\alpha n}$ and $g'(u)=0$ into  \eqref{bcup1} and into \eqref{bcup2}. We obtain
%%%%%%%%%%%%%%%%%%%%%%%%%%%%%%%%%%%%%%%%
\begin{eqnarray}
\label{uprtball1}\Delta u'&=&0\qquad\hbox{in}\:B_R\\
\label{uprtball2}\partial_{\nu}u'+\alpha u'&=&\left(\frac{1+\alpha\:R}{n}\right)(v
\cdot\nu)\qquad\hbox{in}\:\partial B_R.
\end{eqnarray}
2. Similarly we get for the eigenvalue problem in $B_R$
\begin{eqnarray}
\label{uprewball1}&&\Delta u'+\lambda(0)u'+\dot{\lambda}(0)u'=0\qquad\hbox{in}\:B_R\\
\label{uprewball2}&&\partial_{\nu}u'+\alpha u'=\left([1+\alpha\:R-n]\alpha+
\lambda\:R\right)\frac{u(R)}{R}(v\cdot\nu)\quad\hbox{in}\:\partial B_R.
\end{eqnarray}
\end{example}
%%%%%%%%%%%%%%%%%%%%%%%%%%%%%%%%%%%%%%
%%%%%%%%%%%%%%%%%%%%%%%%%%%%%%%%%%%%%
\section{The second domain variation}
\label{secondvar}
%%%%%%%%%%%%%%%%%%%%%%%%%%%%%%%%%%%%%%
%%%%%%%%%%%%%%%%%%%%%%%%%%%%%%%%%%%%%%
The aim of this section is to find a suitable form  of $\ddot{\cal{E}}(0)$ in order to determine its sign. Recall that
$\ddot{\cal{E}}(t)$ is given by \eqref{eddot2} and that consequently
\begin{eqnarray}\label{eddot21}
\ddot{\cal{E}}(0)&=&\int_\Omega \nabla  u \ddot A \nabla  u \:dx +
\alpha \oint_{\p \Omega}  u^2\ddot m\:dS -2\int_\Omega G( u)\ddot J\:dx\\
\nonumber&&\qquad-2\int_{\Omega }|\nabla \dot{\tilde u}|^2\:dx
-2\alpha \oint_{\p \Omega}\dot{\tilde u}^2\:dS +2\int_\Omega g' ( u)\dot
{\tilde u}^2\:dx. 
\end{eqnarray}
In this section we do not assume that $\Omega$ is a critical domain. 
\newline
\newline
The following integrals, which appear in  \eqref{eddot21}, will be expanded 
with respect to $t$.
\begin{eqnarray}
\label{parten1}{\cal{F}}_{1}(t):=\int_{\Omega}\ddot{A}_{ij}(t)\:\partial_{i}
\tilde{u}(t)
\:\partial_{j}\tilde{u}(t)\:dx\\
\label{parten2}{\cal{F}}_{2}(t):=\alpha\oint_{\p\Omega}\tilde{u}^2(t)
\:\ddot{m}(t)\:dS\\
\label{parten3}{\cal{F}}_{3}(t):=-2\int_{\Omega}G(\tilde{u}(t))
\:\ddot{J}(t)\:dx\\
\label{parten4}{\cal{F}}_{4}(t):=-2\int_{\Omega}{A}_{ij}(t)\:\partial_{i}
\dot{\tilde{u}}(t)
\:\partial_{j}\dot{\tilde{u}}(t)\:dx\\
\label{parten5}{\cal{F}}_{5}(t):= -2\alpha \oint_{\p \Omega}\dot{\tilde u}^2(t)m(t)\:dS
\\
\label{parten6}{\cal{F}}_{6}(t):=2\int_\Omega g'(\tilde u)\dot {\tilde u}^2(t)
 J(t)\:dx.
\end{eqnarray}
%%%%%%%%%%%%%%%%%%%%%%%%%%%%%%%%%%%%%%
\subsection{The expression ${\cal{F}}_{1}(0)+{\cal{F}}_{4}(0)$}
\label{subf1f4}
%%%%%%%%%%%%%%%%%%%%%%%%%%%%%%%%%%%%%%
From Lemma \ref{aij3} we have
\begin{eqnarray*}
{\cal{F}}_{1}(0)
&=&
\int_{\Omega}\left((\Div v)^2-D_v:D_v\right)\vert\nabla u\vert^2\:dx\\
&&
+
2\int_{\Omega}\left(\partial_{k}v_i\:\partial_{j}v_k+\partial_{k}v_j\:\partial_{i}v_k\right)\partial_{i}u\:\partial_{j}u\:dx\\
&&+
2\int_{\Omega}\partial_{k}v_{i}\:\partial_{k}v_j\:\partial_{i}u\:\partial_{j}u\:dx\\
&&
-
2\int_{\Omega}\Div v\left(\partial_{j}v_i+\partial_{i}v_{j}\right)\:\partial_{i}u\:\partial_{j}u\:dx+\mathcal{D},
\end{eqnarray*}
where
\begin{eqnarray*}
\mathcal{D}=-\int_\Omega(\p_iw_j+\p_j w_i)\p_i u\p_ju \:dx +\int_\Omega \Div\:w|\nabla u|^2\:dx.
\end{eqnarray*}
We use notation $(D_v)_{ij}=\partial_{j}v_{i}$. We rewrite ${\cal{F}}_{1}(0)$ as
\begin{eqnarray*}
{\cal{F}}_{1}(0)
&=&
\int_{\Omega}\left((\Div v)^2-D_v:D_v\right)\vert\nabla u\vert^2\:dx
+
4\int_{\Omega}(\nabla u\cdot D_v)\cdot(D_v\nabla u)\:dx\\
&&+
2\int_{\Omega}(D_v\nabla u)\cdot(D_v\nabla u)\:dx
-
4\int_{\Omega}\Div v\:\nabla u\cdot D_v\nabla u\:dx+\mathcal{D}.
\end{eqnarray*}
From Lemma \ref{aij3} and \eqref{utilde2} we also have
\begin{eqnarray*}
{\cal{F}}_{4}(0)&:=&-2\int_{\Omega}\:\partial_{i}\dot{\tilde{u}}(0)
\:\partial_{i}\dot{\tilde{u}}(0)\:dx\\
&=&
-
2\int_{\Omega}\partial_{i}v_{k}\:\partial_{k}u\:\partial_{i}v_{l}\:\partial_{l}u\:dx
-
2\int_{\Omega}v_{k}\:\partial_{k}\partial_{i}u\:v_{l}\:\partial_{l}\partial_{i}u
\:dx\\
&&-2
\int_{\Omega}\vert\nabla u'\vert^2\:dx
-
4\int_{\Omega}\:v_{l}\:\partial_{l}\partial_{i}u\:\partial_{i}v_{k}\:\partial_{k}u\:dx\\
&&
-
4\int_{\Omega}\partial_{i}u'\:\partial_{i}v_{k}\:\partial_{k}u\:dx
-
4\int_{\Omega}v_{k}\:\partial_{k}\partial_{i}u\:\partial_{i}u'\:dx.
\end{eqnarray*}
Moreover, in terms of matrices we have, setting $(D^2u)_{ij}=\p_i\p_ju$,
\begin{eqnarray*}
{\cal{F}}_{4}(0)&:=&-2\int_{\Omega}\:\partial_{i}\dot{\tilde{u}}(0)
\:\partial_{i}\dot{\tilde{u}}(0)\:dx
=
-
2\int_{\Omega}(D_v\nabla u)\cdot(D_v\nabla u)\:dx\\
&&
-
2\int_{\Omega}(D^2u\:v)\cdot(D^2u\:v)\:dx
-2
\int_{\Omega}\vert\nabla u'\vert^2\:dx\\
&&
-
4\int_{\Omega}\:(D^2u\:v)\cdot(D_v\nabla u)\:dx
-
4\int_{\Omega}\nabla u'\cdot (D_v \nabla u)\:dx\\
&&
-
4\int_{\Omega}(D^2u\:v)\cdot\nabla u'\:dx.
\end{eqnarray*}
For the sum ${\cal{F}}_{1}(0)+{\cal{F}}_{4}(0)$ we observe that the integral $2\int_
{\Omega}(D_v\nabla u)\cdot(D_v\nabla u)\:dx$ cancels:
\begin{eqnarray*}
{\cal{F}}_{1}(0)+{\cal{F}}_{4}(0)
&=&
\int_{\Omega}\left((\Div v)^2-D_v:D_v\right)\vert\nabla u\vert^2\:dx\\
&&
+
4\int_{\Omega}(\nabla u\cdot D_v)\cdot(D_v\nabla u)\:dx\\
&&
-
4\int_{\Omega}\Div v\:\nabla u\cdot D_v\nabla u\:dx
-
2\int_{\Omega}(D^2u\:v)\cdot(D^2u\:v)\:dx\\
&&-
4\int_{\Omega}\:(D^2u\:v)\cdot(D_v\nabla u)\:dx
-
2\int_{\Omega}\vert\nabla u'\vert^2\:dx\\
&&
-
4\int_{\Omega}\nabla u'\cdot (D_v \nabla u)\:dx
-
4\int_{\Omega}(D^2u\:v)\cdot\nabla u'\:dx+\mathcal{D}.
\end{eqnarray*}
Observe that the last two integrals can be written as
\begin{eqnarray*}
&&-4\int_{\Omega}\nabla u'\cdot (D_v \nabla u)\:dx
-
4\int_{\Omega}(D^2u\:v)\cdot\nabla u'\:dx\\
&&\qquad\qquad=
4\int_{\Omega}v\cdot\nabla u\:\Delta u'\:dx
-
4\oint_{\p\Omega}v\cdot\nabla u\:\partial_{\nu}u'\:dS.
\end{eqnarray*}
We will show that ${\cal{F}}_{1}(0)+{\cal{F}}_{4}(0)$ can be written as a sum of 
boundary integrals and two domain integrals involving the Laplace operator. The 
computations are done in three steps.
\newline
\newline
{\bf{Step 1}} We observe that
\begin{eqnarray*}
I&:=&-4\int_{\Omega}\Div v\:\nabla u\cdot(D_v\nabla u)\:dx
-
4\int_{\Omega}(v\cdot D^2u)\cdot(D_v\nabla u)\:dx\\
&=&
-
4\int_{\Omega}\partial_{j}\left(v_{j}\:\partial_{i}u\right)\:\partial_{i}v_{k}\:\partial_{k}u
\:dx\\
&=&
4\int_{\Omega}v_{j}\:\partial_{i}u\:\partial_{i}\partial_{j}v_k\:\partial_{k}u\:dx
+
4\int_{\Omega}v_{j}\:\partial_{i}u\:\partial_{i}v_k\:\partial_{j}\partial_{k}u\:dx\\
&&-
4\oint_{\p\Omega}(v\cdot\nu)\: \nabla u\cdot(D_v\nabla u)\:dS.
\end{eqnarray*}
We integrate again  the integral $4\int_{\Omega}v_{j}\:\partial_{i}u\:
\partial_{i}\partial_{j}v_k\:\partial_{k}u\:dx$ by parts. This gives a term with $\Delta u$:
\begin{eqnarray*}
I&:=&
-4\int_{\Omega}\Delta u\:v\cdot(D_v\nabla u)\:dx
-
4\int_{\Omega}(\nabla u\cdot D_v)\cdot(D_v\nabla u)\:dx\\
&&-4\int_{\Omega}(v\cdot D_v)\cdot(D^2u\nabla u)\:dx
+
4\int_{\Omega}(\nabla u\cdot D_v)\cdot(D^2 u\:v)\:dx\\
&&+
4\oint_{\p\Omega}\partial_{\nu}u\:v\cdot(D_v\nabla u)\:dS
-
4\oint_{\p\Omega}(v\cdot\nu)\: \nabla u\cdot(D_v\nabla u)\:dS.
\end{eqnarray*}
Then
\begin{eqnarray*}
{\cal{F}}_{1}(0)+{\cal{F}}_{4}(0)
&=&
\int_{\Omega}\left((\Div v)^2-D_v:D_v\right)\vert\nabla u\vert^2\:dx
-
4\int_{\Omega}\Delta u\:v\cdot(D_v\nabla u)\:dx
\\
&&
-
4\int_{\Omega}(v\cdot D_v)\cdot(D^2u\nabla u)\:dx
+
4\int_{\Omega}(\nabla u\cdot D_v)\cdot(D^2 u\:v)\:dx\\
&&
-
2\int_{\Omega}(D^2u\: v)\cdot(D^2u\: v)\:dx
+
4\oint_{\p\Omega}\partial_{\nu}u\:v\cdot(D_v\nabla u)\:dS\\
&&
-
4\oint_{\p\Omega}(v\cdot\nu)\: \nabla u\cdot(D_v\nabla u)\:dS
-
2\int_{\Omega}\vert\nabla u'\vert^2\:dx\\
&&+
4\int_{\Omega}v\cdot\nabla u\:\Delta u'\:dx
-
4\oint_{\p\Omega}v\cdot\nabla u\:\partial_{\nu}u'\:dS+\mathcal{D}.
\end{eqnarray*}
{\bf{Step 2}} Again by partial integration we get
\begin{eqnarray*}
&&-
2\int_{\Omega}(v\cdot D_v)\cdot(D^2u\nabla u)\:dx
-
2\int_{\Omega}(D^2u\: v)\cdot(D^2u\: v)\:dx\\
&&\qquad=
2\int_{\Omega}\Div v\:v\cdot(D^2u\nabla u)\:dx
+
2\int_{\Omega}v_i\:v_j\:\partial_{k}u\:\partial_{i}\partial_{j}\partial_{k}u\:dx\\
&&\qquad\quad-
2\oint_{\p\Omega}(v\cdot\nu)\:v\cdot(D^2u \nabla u)\:dS.
\end{eqnarray*}
Moreover,
\begin{eqnarray*}
&&4\int_{\Omega}(\nabla u\cdot D_v)\cdot(D^2 u\:v)\:dx
=
-
2\int_{\Omega}\Delta u\:v\cdot(D^2u\:v)\:dx\\
\qquad&&
-
2\int_{\Omega}v_i\:v_j\:\partial_{k}u\:\partial_{i}\partial_{j}\partial_{k}u\:dx
+
2\oint_{\p\Omega}\partial_{\nu}u\:v\cdot(D^2u\:v)\:dS.
\end{eqnarray*}
Thus,
\begin{eqnarray*}
{\cal{F}}_{1}(0)+{\cal{F}}_{4}(0)
&=&
\int_{\Omega}\left((\Div v)^2-D_v:D_v\right)\vert\nabla u\vert^2\:dx
-
4\int_{\Omega}\Delta u\:v\cdot(D_v\nabla u)\:dx
\\
&&
-
2\int_{\Omega}(v\cdot D_v)\cdot(D^2u\nabla u)\:dx
+
2\int_{\Omega}\Div v\:v\cdot(D^2u\nabla u)\:dx\\
&&
-
2\int_{\Omega}\Delta u\: v\cdot(D^2u\: v)\:dx
-
2\oint_{\p\Omega}(v\cdot\nu)\:v\cdot(D^2u\:\nabla u)\:dS\\
&&
+
2\oint_{\p\Omega}\partial_{\nu}u\:v\cdot(D^2u\:v)\:dS
+
4\oint_{\p\Omega}\partial_{\nu}u\:v\cdot(D_v\nabla u)\:dS\\
&&
-
4\oint_{\p\Omega}(v\cdot\nu)\: \nabla u\cdot(D_v\nabla u)\:dS
-
2\int_{\Omega}\vert\nabla u'\vert^2\:dx\\
&&
+
4\int_{\Omega}v\cdot\nabla u\:\Delta u'\:dx
-
4\oint_{\p\Omega}v\cdot\nabla u\:\partial_{\nu}u'\:dS+\mathcal{D}.
\end{eqnarray*}
{\bf{Step 3}} Finally we note that
\begin{eqnarray*}
&&\Div\left([v\:\Div v-v\cdot D_v]\vert\nabla u\vert^2\right)\\
&&\qquad=
\left((\Div v)^2-D_v:D_v\right)\vert\nabla u\vert^2
+
2(D^2u\nabla u)\cdot(v\:\Div v-v\cdot D_v).
\end{eqnarray*}
Straightforward partial integration implies
\begin{align}\label{eq:D}
\mathcal{D}=& -2\oint_{\p\Omega} w_i\p_iu\p_\nu u\:dS +\oint_{\p\Omega} (w\cdot\nu) |\nabla u|^2\:dS\\
\nonumber &-2\oint_{\p \Omega} (w\cdot \nu)G(u)\:dS +2\int_\Omega G(u)\rm{div}\;w\:dx.
\end{align}
In summary we have proved
%%%%%%%%%%%%%%%%%%%%%%%%%%%%%%%%%%%%%%
\begin{proposition}\label{propf1f4}
A formal computation without any further assumption on $v$ yields
\begin{eqnarray*}
{\cal{F}}_{1}(0)+{\cal{F}}_{4}(0)
&=&\int_{\Omega}\Div\left([v\:\Div v-v\cdot D_v]\vert\nabla u\vert^2\right)\:dx
-
4\int_{\Omega}\Delta u\:v\cdot(D_v\nabla u)\:dx\\
&&
-
2\int_{\Omega}\Delta u\:(D^2u\:v)\cdot v\:dx
+
4\oint_{\p\Omega}\partial_{\nu}u\:v\cdot(D_v\nabla u)\:dS\\
&&
-
4\oint_{\p\Omega}(v\cdot\nu)\:\nabla u\cdot(D_v\nabla u)\:dS
+
2\oint_{\p\Omega}(D^2 u\:v)\cdot v\:\partial_{\nu}u\:dS\\
&&-
2\oint_{\p\Omega}(v\cdot\nu)\: v\cdot(D^2u\nabla u)\:dS
+
4\int_{\Omega}v\cdot\nabla u\:\Delta u'\:dx\\
&&
-
4\oint_{\p\Omega}v\cdot\nabla u\:\partial_{\nu}u'\:dS
-
2\int_{\Omega}\vert\nabla u'\vert^2\:dx+\mathcal{D}.
\end{eqnarray*}
\end{proposition}
%%%%%%%%%%%%%%%%%%%%%%%%%%%%%%%%%%%%%%
\subsection{The expression ${\cal{F}}_{3}(0)+{\cal{F}}_{6}(0)$}
\label{subf3f6}
%%%%%%%%%%%%%%%%%%%%%%%%%%%%%%%%%%%%%%
From \eqref{parten3}, \eqref{jform1} and \eqref{utilde1} we have
\begin{eqnarray*}
{\cal{F}}_{3}(0)&=&-2\int_{\Omega}G(\tilde{u}(0))\:\ddot{J}(0)\:dx\\
&=&
-2\int_{\Omega}G(u(x))\:\left((\Div v)^2- D_v:D_v+\Div\:w\right)\:dx.
\end{eqnarray*}
We again use the fact that
\begin{eqnarray*}
(\Div v)^2- D_v:D_v=\Div\left(v\:\Div v-v\cdot D_v\right).
\end{eqnarray*}
We get
\begin{eqnarray*}
{\cal{F}}_{3}(0)&=&2\int_{\Omega}g(u(x))\left(v\cdot\nabla u\:\Div v-v\cdot (D_v\nabla u)\right)\:dx\\
&&-
2\oint_{\p\Omega}G(u(x))\left((v\cdot\nu)\:\Div v-v\cdot (D_v\nu)\right)\:dS\\
&&-
2\int_\Omega G(u)\Div\:w\:dx.
\end{eqnarray*}
From \eqref{parten6}, \eqref{jform3}  and \eqref{utilde1} we have
\begin{eqnarray*}
{\cal{F}}_{6}(0)&=&2\int_\Omega g'(\tilde{u}(0))\:\dot{\tilde u}^2(0)\:J(0)\:dx
=
2\int_\Omega g'(u(x))\left(v\cdot\nabla u+u'\right)^2\:dx\\
&=&
2\int_\Omega g'(u(x))\left(v(x)\cdot\nabla u(x)\right)^2\:dx
+
2\int_\Omega g'(u(x))\:u'^2(x)\:dx
\\
&&+
4\int_\Omega g'(u(x))\:v(x)\cdot\nabla u(x)\:u'(x)\:dx.
\end{eqnarray*}
We note that
\begin{eqnarray*}
&&2\int_\Omega g'(u(x))\left(v(x)\cdot\nabla u(x)\right)^2\:dx
=
2\int_\Omega v\cdot\nabla g(u)\:v\cdot\nabla u\:dx\\
&&\qquad=
-
2\int_\Omega g(u(x))\Div\:v(x)\:v(x)\cdot\nabla u\:dx\\
&&\qquad\quad
-
2\int_\Omega g(u(x))\:v(x)\cdot(D_v\nabla u(x))\:dx\\
&&\qquad\quad
-
2\int_\Omega g(u(x))\:v(x)\cdot(D^2 u(x)v(x))\:dx\\
&&\qquad\quad
+
2\oint_{\p\Omega} g(u(x))\:v(x)\cdot\nu\:v\cdot\nabla u(x)\:dS.
\end{eqnarray*}
From this we easily deduce the following proposition.
%%%%%%%%%%%%%%%%%%%%%%%%%%%%%%%%%%%%%%
\begin{proposition}\label{propf3f6}
A formal computation without any further assumption on $v$ yields
\begin{eqnarray*}
{\cal{F}}_{3}(0)+{\cal{F}}_{6}(0)
&=&
-4\int_{\Omega}g(u(x))\:v(x)\cdot(D_v\nabla u(x))\:dx\\
&&-
2\int_{\Omega}g(u(x))\:v(x)\cdot(D^2 u(x)\:v(x))\:dx\\
&&-
2\oint_{\p\Omega}G(u)\:\left(v(x)\cdot\nu\:\Div v-v(x)\cdot(D_v\nu)\right)
\:dS\\
&&+
2\oint_{\p\Omega} g(u(x))\:v(x)\cdot\nu\:v\cdot\nabla u(x)\:dS\\
&&
+
2\int_{\Omega}g'(u(x))\:u'^{2}(x)\:dx-
2\int_\Omega G(u)\Div\:w\:dx\\
&&+
4\int_{\Omega} g'(u(x))\:v(x)\cdot\nabla u(x)\:u'(x)\:dx.
\end{eqnarray*}
\end{proposition}
%%%%%%%%%%%%%%%%%%%%%%%%%%%%%%%%%%%%%%
\subsection{The expression ${\cal{F}}_{2}(0)+{\cal{F}}_{5}(0)$}
\label{subf2f5}
%%%%%%%%%%%%%%%%%%%%%%%%%%%%%%%%%%%%%%
From \eqref{parten2} and \eqref{utilde1} we deduce
\begin{eqnarray*}
{\cal{F}}_{2}(0):=\alpha\oint_{\p\Omega}\tilde{u}^2(0)
\:\ddot{m}(0)\:dS
=
\alpha\oint_{\p\Omega}u^2(x)\:\ddot{m}(0)\:dS.
\end{eqnarray*}
We will not use the explicit for of $\ddot{m}(0)$. From \eqref{utilde2} and the fact that $m(0)=1$ we obtain
\begin{eqnarray*}
{\cal{F}}_{5}(0)&:=& -2\alpha\oint_{\p \Omega}\dot{\tilde u}^2(0)m(0)\:dS\\
&=&
-
2\alpha\oint_{\p \Omega}(v\cdot\nabla u)^2\:dS
-
4\alpha\oint_{\p \Omega}v\cdot\nabla u\:u'\:dS
-
2\alpha\oint_{\p \Omega}u'^{2}\:dS.
\end{eqnarray*}
Thus,
\begin{eqnarray}\label{secsum}
{\cal{F}}_{2}(0)+{\cal{F}}_{5}(0)
&=&
-
2\alpha\oint_{\p \Omega}(v\cdot\nabla u)^2\:dS
-
4\alpha\oint_{\p \Omega}v\cdot\nabla u\:u'\:dS\\
\nonumber&&
-
2\alpha\oint_{\p \Omega}u'^{2}\:dS
+
\alpha\oint_{\p\Omega}u^2(x)\:\ddot{m}(0)\:dS.
\end{eqnarray}
%%%%%%%%%%%%%%%%%%%%%%%%%%%%%%%%%%%%%%
\subsection{Main result}
\label{mr}
We add up all these contributions. We then arrive at our final result.
%%%%%%%%%%%%%%%%%%%%%%%%%%%%%%%%%%%%%%
\begin{theorem}\label{ddE}
Assume that $\Delta u+ g(u)=0$ in $\Omega$ and $\p_\nu u + \alpha u=0$ on $\p \Omega$.  Let  $u'$ satisfy \eqref{bcup1} and \eqref{bcup2}. Then the second variation $\ddot{{\cal{E}}}(0)$ can be expressed in the form
\begin{eqnarray}\label{ddef}
\nonumber\ddot{{\cal{E}}}(0)&=&
+
\oint_{\p\Omega}[(v\cdot\nu)\:\Div v-v\cdot(D_v\nu)+w\cdot\nu]\left(\vert\nabla u\vert^2-
2G(u)\right)\:dS\\
\nonumber&&+
4\oint_{\p\Omega}\left(\partial_{\nu}u\:v\cdot(D_v\nabla u)-(v\cdot\nu)\:\nabla u
\cdot(D_v\nabla u)\right)\:dS\\
&&+
2\oint_{\p\Omega}\left(\partial_{\nu}u\:v\cdot(D^2u\: v)-(v\cdot\nu)\:v
\cdot(D^2u \nabla u)\right)\:dS\\
\nonumber&&+
2\oint_{\p\Omega}g(u)\:(v\cdot\nu)\:v\cdot\nabla u\:dS
-
4\oint_{\p\Omega}v\cdot\nabla u\left(\partial_{\nu}u'+\alpha u'\right)\:dS\\
\nonumber&&
-
2\alpha\oint_{\p\Omega}(v\cdot\nabla u)^2\:dS
-
2\oint_{\p\Omega}w\cdot\nabla u\:\partial_{\nu}u\:dS 
\\
\nonumber&&
+
\alpha\oint_{\p\Omega}u^2(x)\:\ddot{m}(0)\:dS
-
2Q_g(u'),
\end{eqnarray}
where
\begin{eqnarray}\label{quadratic}
Q_g(u'):= \int_\Omega |\nabla u'|^2\:dx -\int_\Omega g'(u)u'^2\:dx +\alpha \oint_{\p \Omega} u'^2\:dS
\end{eqnarray}
is a quadratic form in $u'$.
\end{theorem}
%%%%%%%%%%%%%%%%%%%%%%%%%%%%%%%%%%%%%
This formula is very general because no volume constraint is used. It could for instance be used to study problems with a prescribed surface area.
%%%%%%%%%%%%%%%%%%%%%%%%%%%%%%%%%%%%%
\section{Applications to nearly spherical domains}
\label{spherical}
%%%%%%%%%%%%%%%%%%%%%%%%%%%%%%%%%%%%%
\subsection{Second variation}
\label{subsecvar}
%%%%%%%%%%%%%%%%%%%%%%%%%%%%%%%%%%%%%
We evaluate \eqref{ddef} if 
$\Omega=B_R$, $u=u(\vert x\vert)$ and when the domain perturbations preserve the volume  and satisfy \eqref{volcomp1} and \eqref{volume2}. Then
\begin{eqnarray*}
\partial_{i}\partial_{j}u(\vert x\vert)\vert_{\partial B_R}
=
\partial_{\nu}^2u(R)\:\nu_{i}\:\nu_{j}+\frac{1}{R}\partial_{\nu}u(R)\left(\delta_{ij}-\nu_{i}\:\nu_{j}\right).
\end{eqnarray*}
Since $\vert\nabla u\vert^2-2G(u)=const.$ on $\partial B_R$  the contribution in the first integral of \eqref{ddef} vanishes by \eqref{volcomp1}. We keep in mind the Robin boundary condition for $u$.  We get
\begin{eqnarray}\label{ddefrad}
&&\ddot{{\cal{E}}}(0)=
-
2Q_g(u')+
\alpha\:u^2(R)\oint_{\p B_R}\:\ddot{m}(0)\:dS\\
\nonumber&&
+
4\alpha^2\:u^2(R)\oint_{\p B_R}v^\tau\cdot D_v\:\nu\:dS
+ 
\frac{2\alpha^2}{R}\:u^2(R)\oint_{\p B_R}(v^\tau)^2\:dS\\
\nonumber&&-
2\alpha\:u(R)\oint_{\p B_R}g(u)\:(v\cdot\nu)^2\:dS-
2u_r(R)^2\oint_{\p B_R}w\cdot\nu\:dS\\
\nonumber&&
+
4\alpha\:u(R)\oint_{\p B_R}(v\cdot\nu)\left(\partial_{\nu}u'+\alpha u'\right)\:dS
-
2\alpha^3\:u^2(R)\oint_{\p B_R}(v\cdot\nu)^2\:dS.
\end{eqnarray}
We need the following technical lemma for $v^\tau$.
%%%%%%%%%%%%%%%%%%%%%%%%%%%%%%%%%%%%%%%%
\begin{lemma}\label{techrad} For volume preserving perturbations there holds
\begin{eqnarray*}
&&\frac{2\alpha^2}{R}\:u^2(R)\oint_{\p B_R}(v^{\tau})^2\:dS
=
-
4\alpha^2\:u^2(R)\oint_{\p B_R}v^\tau \cdot D_{v}\nu\:dS\\
&&
\qquad+
\frac{2\alpha^2(n-1)}{R}\:u^2(R)\oint_{\p B_R}\:(v\cdot \nu)^2\:dS
+
2\alpha^2u(R)^2\oint_{\p B_R}w\cdot\nu\:dS.
\end{eqnarray*}
\end{lemma}
%%%%%%%%%%%%%%%%%%%%%%%%%%%%%%%%%%%%%%%%
At first observe that
\begin{eqnarray*}
\frac{1}{R}\oint_{\p B_R}(v^\tau)^2\:dS
=
\oint_{\p B_R}v\cdot D_{\nu}v\:dS.
\end{eqnarray*}
On the other hand, 
\begin{eqnarray*}
\oint_{\p B_R}v\cdot D_{\nu}v\:dS
&=&
\oint_{\p B_R}v^\tau_i(\nabla^\tau_i\nu_k)v_k\:dS\\
&=&
-\oint_{\p B_R}v^\tau_i(\nabla^\tau_iv_k)\nu_k\:dS
-
\oint_{\p B_R}\Div_{\partial B_R}v^\tau\:(v\cdot \nu)\:dS\\
&=&
-
\oint_{\p B_R}v\cdot D_{v}\nu\:dS
-
\oint_{\p B_R}\Div v\:(v\cdot \nu)\:dS\\
&& 
+
2\oint_{\p B_R}(v\cdot\nu)\:\nu\cdot D_v\nu\:dS
+
\oint_{\p B_R}\Div_{\partial B_R}\nu\:(v\cdot \nu)^2\:dS.
\end{eqnarray*}
Next we use \eqref{volcomp1}. Then
\begin{eqnarray*}
\oint_{\p B_R}v\cdot D_{\nu}v\:dS
&=&
-
2\oint_{\p B_R}v^\tau\cdot D_{v}\nu\:dS
+
\frac{n-1}{R}\oint_{\p B_R}\:(v\cdot \nu)^2\:dS\\
&&+
\oint_{\p\Omega}w\cdot\nu\:dS.
\end{eqnarray*}
This proves the claim.
\hfill{$\rlap{$\sqcap$}\sqcup$}
\newline\newline
%%%%%%%%%%%%%%%%%%%%%%%%%%%%%%%%%%%%%%%%
This lemma, together with the Robin condition $u_r(R)+\alpha u(R)=0$, implies that \eqref{ddefrad} can be written as
\begin{eqnarray}\label{ddefrad1}
&&\ddot{{\cal{E}}}(0)=
-
2Q_g(u')\\
\nonumber&&
+
\alpha\:u^2(R)\oint_{\p B_R}\:\ddot{m}(0)\:dS
+
4\alpha\:u(R)\oint_{\p B_R}(v\cdot\nu)\left(\partial_{\nu}u'+\alpha u'\right)\:dS\\
\nonumber&&
-
2\alpha\:u(R)\left(g(u)-\frac{\alpha(n-1)}{R}\:u(R) -\alpha^2 u(R)\right)\oint_{\p B_R}\:(v\cdot \nu)^2\:dS.
\end{eqnarray}
We rewrite \eqref{bcup2} for the radial situation. Recall that
\begin{eqnarray*}
\partial_{\nu}u'+\alpha u'&=& -\p_\nu( v\cdot\nabla u)+ 
\nabla^\tau u\cdot D_v\nu+\nu \cdot D_v\nabla u
-\alpha \:v\cdot\nabla u\qquad\hbox{in}\:\partial\Omega.
\end{eqnarray*}
Thus, in the radial case we get
\begin{eqnarray}
\label{genform1}\Delta u'+g'(u)u'&=&0\quad\hbox{in}\:B_R\\
\label{genform2}\partial_{\nu}u'+\alpha u'&=&k_g(u(R))(v\cdot\nu)\quad\hbox{in}\:\partial B_R,
\end{eqnarray}
with
\begin{eqnarray}\label{defk}
k_g(u(R)):= g(u(R))-\frac{\alpha(n-1)}{R}u(R)+\alpha^2u(R).
\end{eqnarray}
We can insert \eqref{genform2} into \eqref{ddefrad1} and obtain $\ddot{{\cal{E}}}(0)$ as a quadratic functional in $u'$ alone.
\begin{eqnarray}\label{ddefrad2}
\ddot{{\cal{E}}}(0)&=&
-
2Q_g(u')+
\alpha\:u^2(R)\oint_{\p B_R}\:\ddot{m}(0)\:dS\\
\nonumber&&+
\frac{2\alpha\:u(R)}{k_g(u(R))}\oint_{\p B_R}\left(\partial_{\nu}u'+\alpha u'\right)^2\:dS.
\end{eqnarray}
Further simplification is possible if we use \eqref{bcup1} for $u'$. We multiply this equation with $u'$ and integrate over $B_R$. This leads to
\begin{eqnarray*}
Q_g(u')= \oint_{\p B_R}(\partial_{\nu}u'+\alpha u')u'\:dS.
\end{eqnarray*}
%%%%%%%%%%%%%%%%%%%%%%%%LLLLLLLLLLLLLLLL
\begin{lemma} For every volume preserving perturbation of the ball and for radially symmetric solutions  $u$ we have
\begin{eqnarray}\label{Eradial}
\ddot{{\cal{E}}}(0)&=&-2 \oint_{\p B_R}(\partial_{\nu}u'+\alpha u')u'\:dS
+
\alpha\:u^2(R)\oint_{\p B_R}\:\ddot{m}(0)\:dS \\
\nonumber&&\qquad
+
\frac{2\alpha\:u(R)}{k_g(u(R))}\oint_{\p B_R}\left(\partial_{\nu}u'+\alpha u'\right)^2\:dS.
\end{eqnarray}
\end{lemma}
\begin{remark} The second variation is independent of $v^\tau$ and $w$. We can therefore restrict ourselves to Hadamard perturbations $y=x+tN\nu +O(t^2)$.
\end{remark}
Consider the case where $(v\cdot \nu)=0$ on $\p B_R$. Then by \eqref{genform1}, \eqref{genform2} and \eqref{quadratic} we have $Q_g(u')=0$. Moreover, by  \eqref{isodefrad1}, Lemma \ref{Hadamard}
\eqref{ddefrad2} it follows that $\ddot{{\cal{E}}}(0)=0$. Consequently perturbations, which preserve the volume and for which $(v\cdot \nu)=0$, lie in the kernel of $\ddot{{\cal{E}}}(0)$.
%%%%%%%%%%%%%%%%%%%%%%%%%%%%%%%%%%%%%%%%
\subsection{Discussion of the sign of $\ddot{{\cal{E}}}(0)$ in the radial case }
\label{subsign}
%%%%%%%%%%%%%%%%%%%%%%%%%%%%%%%%%%%%%%%%
\subsubsection{General strategy}
\label{subsubgs}
%%%%%%%%%%%%%%%%%%%%%%%%%%%%%%%%%%%%%%%%
We recall \eqref{ddefrad2} and \eqref{genform1}. Then
\begin{align}\label{ddE1}
\ddot{{\cal{E}}}(0)=\alpha u^2(R)\ddot{\cal{S}}(0) +{\cal{F}},
\end{align}
where
\begin{align*}
{\cal{F}}:=-2Q_g(u')+2\alpha\:u(R)k_g(u(R)) \oint_{\p B_R}(v\cdot \nu)^2\:dS.
\end{align*}
By Lemma \ref{Hadamard} we have $\ddot{\cal{S}}(0)>0$. In order to estimate ${\cal{F}}$  we consider the following Steklov eigenvalue problem
\begin{align}\label{steklovg}
\Delta \phi +g'(u)\phi=0 \tx{in} B_R,\\
\nonumber \partial_{\nu} \phi +\alpha \phi =\mu \phi \tx{on} \p B_R.
\end{align}
If $g'(u)$ is bounded, then exists an infinite number of eigenvalues
\begin{align*}
\mu_1<\mu_2\leq \mu_3\leq ... \lim_{i\to \infty} \mu_i=\infty
\end{align*}
and a complete system of eigenfunctions $\{\phi_i\}_{i\geq 1}$. We use $\phi_j$ as a test function in \eqref{steklovg}. We find
\begin{align*}
\int_{B_R} [-\nabla \phi_i\cdot \nabla \phi_j + g'(u) \phi_i\phi_j]\:dx -\alpha \oint_{\p B_R} \phi_i\phi_j\:dS+\mu_i\oint_{\p B_R} \phi_i\phi_j\:dS=0.
\end{align*}
If we interchange $i$ and $j$, then we see immediately that the system of eigenfunctions $\{\phi_i\}_i$ can be chosen such that
\begin{eqnarray}\label{orthogonalg}
&\oint_{\p B_R}\phi_i \phi_j\:dS=\delta_{ij}
\end{eqnarray}
and
\begin{eqnarray*}
q(\phi_i,\phi_j)&:=&\int_{B_R} \nabla \phi_i\cdot \nabla \phi_j \:dx -\int_{B_R}g'(u)\phi_i\phi_j \:dx\\
&& +\alpha \oint_{\p B_R} \phi_i\phi_j\:dS = \mu_i\delta_{ij}.
\end{eqnarray*}
We write
\begin{align*}
u'(x)=\sum_{i=1}^\infty c_i \phi_i \quad\tx{and} \quad(v\cdot \nu)=\sum_{i=1}^\infty b_i\phi_i.
\end{align*}
Note that the first eigenfunction $\phi_1$ is radially symmetric and does not change.
The condition 
\begin{eqnarray*}
0=\oint_{\p B_R}(v\cdot \nu)\:dS=\oint_{\p B_R} \phi_1 (v\cdot \nu)\:dS
\end{eqnarray*}
implies that $b_1=0$. It gives a condition on $c_1\:\mu_1$ if we take into account \eqref{genform1} - \eqref{genform2}:
\begin{eqnarray*}
0&=&\oint_{\p B_R}\left(\partial_{\nu}u'+\alpha u'\right)\:\phi_1\:dS\\
&=&
\int_{B_R}\Delta u'\:\phi_1\:dx+\int_{B_R}\nabla u'\nabla\phi_1\:dx
+
\alpha\oint_{\p B_R}u'\:\phi_1\:dS\\
&=&
\oint_{\p B_R}u'\:\partial_{\nu}\phi_1\:dS+\alpha\oint_{\p B_R}u'\:\phi_1\:dS\\
&=&
\mu_1\oint_{\p B_R}u' \:\phi_1\:dS=c_1\:\mu_1.
\end{eqnarray*}
Thus, $c_1\:\mu_1=0$.
The coefficients $b_i$ for  $i\geq 2$ are determined from the boundary value 
problem \eqref{bcup1}, \eqref{bcup2}. In fact,
\begin{eqnarray*}
b_i= \frac{c_i\:\mu_i}{k_g(u(R))}\qquad\hbox{for} \quad i=2,3\dots .
\end{eqnarray*}
By means of the orthonormality conditions of the eigenfunctions we find
\begin{align*}
Q_g(u')&=q(u',u')= \sum_{i=1}^\infty c_i^2q(\phi_i,\phi_i)=\sum_{i=2}^\infty c_i^2\mu_i.
\end{align*}
We insert this into \eqref{ddE1}. We find 
\begin{align}\label{F}
{\cal{F}}= 2\sum_2^\infty c_i^2\:\mu_i^2\left[\frac{\alpha u(R)}{k_g(u(R))}-\frac{1}{\mu_i}\right],
\end{align}
where $k_g$ is defined in \eqref{defk}.
Let $\mu_p=\min\{\mu_i:\mu_i>0\}$ be the smallest  positive eigenvalue. Then
\begin{eqnarray*}
{\cal{F}}&\geq& 2\sum_2^\infty c_i^2\:\mu_i^2\left[\frac{\alpha u(R)}{k_g(u(R))}-\frac{1}{\mu_p}\right]\\
&=&
2k_g^2(u(R))\left[\frac{\alpha u(R)}{k_g(u(R))}-\frac{1}{\mu_p}\right]\oint_{\p B_R}(v\cdot \nu)^2\:dS.
\end{eqnarray*}
The expression ${\cal{F}}$ vanishes if 
\begin{itemize}
\item $u'=0$
\item $\mu_i= \frac{\alpha u(R)}{k_g(u(R))}$ and $(\nu\cdot v)= d_i \phi_i$ for some $i$.
\end{itemize}
The first case occurs only in the case of translations. This, together with Lemma \ref{Hadamard}, implies
%%%%%%%%%%%%%%%%%%%%%%%%%%%%%%%%%%%%%%%%%%LLLLLLLLL
\begin{lemma} The kernel of $\ddot{\cal{E}}(0)$ consists only of first order translations 
\begin{eqnarray*}(\nu\cdot v)=a_ix_i.
\end{eqnarray*}
\end{lemma}
%%%%%%%%%%%%%%%%%%%%%%%%%%%%%%%%%%%%llllllllll
In order to get an estimate for $\ddot{\cal{E}}(0)$ in terms of $v$, we impose the "barycenter" condition
\begin{eqnarray}\label{ortho1}
\oint_{\p B_R}x\:(v(x)\cdot\nu(x))\:dS=0.
\end{eqnarray}
By \eqref{isodefrad1} and \eqref{Stek} it then follows that
\begin{align*}
\oint_{\p B_R} |\nabla^\tau  N|^2\:dS\geq \frac{2n}{R^2}\oint_{\p B_R} N^2\:dS,
\end{align*}
and thus $\ddot{\cal{S}}(0) \geq \frac{n+1}{R^2} \oint_{\p B_R} N^2\:dS.$ Observe that $b_i=0$ for $i=1,\dots, n$. Hence the estimate given above can be improved if we replace $\mu_p$ by $\mu_{p'}=\min \{ \mu_k>0, k>n\}$. This, together with the estimate for ${\cal{F}}$ given above, implies
\begin{align}\label{estimateF}
\ddot{\cal{E}}(0)\geq \left\{ \alpha u^2(R)\frac{n+1}{R^2} +2k_g(u(R))\alpha u(R) -\frac{2k_g^2(u(R))}{\mu_{p'}}\right\}\oint_{\p B_R}(v\cdot \nu)^2\:dS.
\end{align}
In summary we have
%%%%%%%%%%%%%%%%%%%%%%%%%%%%%%%%%%%%%%%%%
\begin{theorem} \label{second}The second variation of $ {\cal{E}}$  for volume preserving perturbations of the first and second order is of the form
\begin{align*}
\ddot{\cal{E}}(0)=\alpha u^2(R)\ddot{\cal{S}}(0) +{\cal{F}}.
\end{align*}
(i) If $\alpha >0$, then it is bounded from below by
\begin{align*}
\alpha u^2(R)\ddot{\cal{S}}(0) + 2k_g^2(u(R))\left[\frac{\alpha u(R)}{k_g(u(R))}-\frac{1}{\mu_p}\right]\oint_{\p B_R}(v\cdot \nu)^2\:dS.
\end{align*}
(ii) Under the additional assumption \eqref{ortho1} we have 
\begin{align*}
\ddot{\cal{E}}(0)\geq \left\{ \alpha u^2(R)\frac{n+1}{R^2} +2k_g(u(R))\alpha u(R) -\frac{2k_g^2(u(R))}{\mu_{p'}}\right\}\oint_{\p B_R}(v\cdot \nu)^2\:dS
\end{align*}
for $\alpha >0$.
\end{theorem}
%%%%%%%%%%%%%%%%%%%%%%%%%%%%%%%%%%%%%%%%
We will apply this theorem to the torsion problem and the eigenvalue problem.
%%%%%%%%%%%%%%%%%%%%%%%%%%%%%%%%%%%%%%%%
\subsubsection{The torsion problem $g=1$}
%%%%%%%%%%%%%%%%%%%%%%%%%%%%%%%%
From \eqref{ballrt} we have $u(R)=\frac{R}{\alpha n}$ and by \eqref{uprtball2} 
\begin{eqnarray*}
k_1= \frac{1+\alpha R}{n} \tx{and}\frac{k_1(u(R))}{\alpha u(R)}=\frac{1+\alpha R}{R}.
\end{eqnarray*} 
The Steklov  problem \eqref{steklovg} is in this case
\begin{align*}
-\Delta \phi =0 \tx{in} B_R, \qquad \p_\nu \phi+\alpha \phi= \mu \phi \tx{on} \p B_R.
\end{align*}
An elementary computation yields $\mu_1-\alpha=0$ and $\mu_k-\alpha=\frac{k-1}{R}$ for $k\geq 2$ and counted without multiplicity. The second eigenvalue $ \mu_2=1/R+\alpha$ has multiplicity $n$ and its  eigenfunctions are
$\frac{x_1}{R}, \hdots, \frac{x_n}{R}$. If $\alpha > 0$, then $\mu_2>0$ and by Theorem \ref{second} (i) 
the inequalities $\ddot{\cal{E}}(0) \geq \ddot{\cal{S}}(0) \geq 0$ hold. Equality holds only for translations. 
\newline
\newline
The estimate can be improved by assuming \eqref{ortho1}. Note that this condition implies in addition to $c_1=0$ also that $c_2=\dots=c_n=0$. Hence we can take $\mu_{p'}= \frac{2}{R}+\alpha$. By Theorem \ref{estimateF} (ii)
\begin{align*}
\ddot{\cal{E}}(0) \geq \left\{\frac{n+1}{\alpha n^2} +2\frac{(1+\alpha R)R}{n^2(2+\alpha R)}\right\}\oint_{\p B_R}(v\cdot \nu)^2\:dS>0.
\end{align*}
Next consider the case where $-\frac{1}{R}<\alpha<0$. Then $\mu_2>0$. Thus, we have $Q(u')>0$. Moreover, ${\cal{F}}<0$ and consequently, by  \eqref{ddE1}, the second variation becomes negative.
\newline
\newline
This property is no longer true if $\alpha \leq -1/R$. In fact, we can always find $c_i$ or equivalently $b_i$, such that ${\cal{F}}>0$ or ${\cal{F}}<0$ and $\ddot{\cal{E}}(0)$ is positive or negative, respectively.  For the torsion problem we have proved the following theorem.
%%%%%%%%%%%%%%%%%%%%%%%%%%%%%%%%%%%%%%%%%%%%%%%
\begin{theorem} (i) Assume $\alpha >0$. Then $\ddot{\cal{E}}(0)>0$ for all volume preserving perturbations.\\
(ii) If $-\frac{1}{R}<\alpha <0$ then $\ddot{\cal{E}}(0)<0$ for all volume preserving perturbations.\\
(iii) If $\alpha\leq -\frac{1}{R}$ then the sign of $\ddot{\cal{E}}(0)$ can change depending on the particular perturbation.
\end{theorem}
This result will be analyzed in \cite{BaWa15} in more detail.
%%%%%%%%%%%%%%%%%%%%%%%%%%%%%%%%%%%%%%%%
\subsubsection{Principal eigenvalue $g(u)=\lambda u$}
%%%%%%%%%%%%%%%%%%%%%%%%%%%%%%%%
In \cite{Ba77} it was shown that for $\alpha_0<\alpha <0$, the ball yields the maximal principal eigenvalue among all nearly spherical domains of the same volume. In this section we therefore restrict ourselves to the case $\alpha >0$. 
\newline
The first eigenfunction is of the form
$u=u(r)= J_{\frac{n-2}{2}}(\sqrt{\lambda}r)r^{-\frac{n-2}{2}}$. Furthermore,
\begin{eqnarray*}
k_{\lambda u}(u(R)):=(\alpha^2R-\alpha(n-1)+\lambda R)\frac{u(R)}{R}. 
\end{eqnarray*}
By \eqref{signk} the term $\alpha^2R-\alpha(n-1)+\lambda R$ is positive. In order to prove that $\ddot\lambda(0)$ is non-negative, we shall use the form \eqref{ddefrad2}.  Under the assumption that $\int_{B_R}u^2\:dx=1$ it follows that
\begin{align}\label{ddlambda}
\ddot{\lambda}(0)=  -2Q_{\lambda u}(u') +2\alpha u(R)k_{\lambda u}\oint_{\p B_R}(v\cdot \nu)^2\:dS+\alpha u^2(R)\ddot{\cal{S}}(0).
\end{align} 
The corresponding Steklov eigenvalue problem is
\begin{eqnarray}
\label{Phieq1}\Delta\phi+\lambda\phi&=&0\qquad\hbox{in}\:B_R\\
\label{Phieq2}\partial_{\nu}\phi+\alpha \phi&=& \mu\phi\qquad\hbox{in}\:\partial B_R.
\end{eqnarray}
Note that $\phi_1= u$ and therefore $\phi_1=$const. on $\p B_R$. Moreover, $\mu_1=0$, therefore $\mu_p=\mu_2$. 
Next we want to check the sign of the expression $\frac{\alpha u(R)}{k_{\lambda u}}-\frac{1}{\mu_2}$ in Theorem \ref{second}. This is equivalent to the sign of
\begin{align*}
L:= \mu_2-\alpha +\frac{n-1}{R} -\frac{\lambda}{\alpha}.
\end{align*}
For this purpose we need the eigenvalues of \eqref{Phieq2}.
The eigenfunctions of $\eqref{Phieq1}, \eqref{Phieq2}$ are of the form
\begin{eqnarray*}
\phi(x)=\sum_{s,i}c_{s,i}\:a_{s}(r)\:Y_{s,i}(\theta),\qquad\theta\in\:S^{n-1}.
\end{eqnarray*}
Here $s\in\N\cup\{0\}$ and $i=1,\hdots,d_s$ for $d_s=(2s+n-2)\:\frac{(s+n-3)!}{s!(n-2)!}\in\N$. The function $Y_{s,i}(\theta)$ denotes the $i$ - th spherical harmonics of order $s$. 
In particular
\begin{eqnarray*}
\Delta^{*}Y_{s,i}+s(s+n-2)Y_{s,i}=0\qquad\hbox{in}\:S^{n-1},
\end{eqnarray*} 
where $\Delta^{*}$ is the Laplace Beltrami operator on the sphere.
As a consequence of this Ansatz we get from \eqref{Phieq1}
\begin{eqnarray*}
a_s(r)=r^{\frac{2-n}{2}}\:J_{s+\frac{n}{2}-1}(\sqrt{\lambda}\:r).
\end{eqnarray*}
The corresponding eigenvalue follows from \eqref{Phieq2}, namely
\begin{eqnarray*}
a_s'(R)=(\mu-\alpha)a_s(R).
\end{eqnarray*}
Since the first eigenfunction does not change sign we have
\begin{align*}
u=r^{-\frac{n-2}{2}}J_{\frac{n-2}{2}}(\sqrt{\lambda}r) .
\end{align*}
It follows from the well-known Bessel identity 
\begin{align*}
(z^{-\nu}J_\nu(z))_z= -r^{-\nu}J_{\nu+1}(z)
\end{align*}
and from $u_r(R)+\alpha u(R)=0$
that
\begin{align}\label{alpha}
\alpha = \sqrt{\lambda}\frac{J_{n/2}(\sqrt{\lambda}R)}{J_{(n-2)/2}(\sqrt{\lambda}R)}.
\end{align}
The eigenfunctions corresponding to $\mu_2$ span the n-dimensional linear space  $(s=1)$
\begin{align*}
\phi(r,\theta)=\sum_{i=1}^nc_{2,i}\:r^{\frac{2-n}{2}}\:J_{\frac{n}{2}}(\sqrt{\lambda}\:r)Y_{1,i}(\theta).
\end{align*}
The boundary condition gives by means of the same identity as before
\begin{align*}
\frac{1}{R}+\alpha-\sqrt{\lambda}\frac{J_{n/2+1}(\sqrt{\lambda}R)}{J_{n/2}(\sqrt{\lambda}R)}=\mu_2.
\end{align*}
If we replace $\alpha$ and $\mu_2$, then we obtain
\begin{align*}
L= \frac{n}{R}-\frac{\sqrt{\lambda}}{J_{\frac{n}{2}}(\sqrt{\lambda}R)}(J_{\frac{n}{2}+1}(\sqrt{\lambda}R)+J_{\frac{n}{2}-1}(\sqrt{\lambda}R)).
\end{align*}
From the identity
\begin{align}\label{Bessel}
nJ_{n/2}(z)=z(J_{n/2+1}(z)+J_{n/2-1}(z)
\end{align}
it follows that $L=0$. 
Consequently for all $v\neq$ const.
\begin{align}\label{f1}
\ddot{\lambda}(0)\geq \alpha u^2(R) \ddot{\cal{S}}(0)>0 .
\end{align}
As for the torsion problem, this inequality can be improved if we impose the barycenter conditions \eqref{volcomp1}.
The positivity of the decond variation is in accordance with Daners-Bossel's inequality \cite{Da06}.
%%%%%%%%%%%%%%%%%%%%%%%%%%%%%%%%%%%%%%%%%%%
%%%%%%%%%%%%%%%%%%%%%%%%%%%%%%%%%%%%%%%%%%%
\section{The ball is optimal}
\label{ballopt}
%%%%%%%%%%%%%%%%%%%%%%%%%%%%%%%%%%%%%%%%%%%
%%%%%%%%%%%%%%%%%%%%%%%%%%%%%%%%%%%%%%%%%%%
By the Taylor expansion
\begin{eqnarray*}
{\cal{E}}(t)={\cal{E}}(0)+t\:\dot{{\cal{E}}}(0)+\frac{t^2}{2}\left(\ddot{{\cal{E}}}(0)+\int_{0}^{\hat{t}}\dddot{{\cal{E}}}(s)\:ds\right)
\end{eqnarray*}
for some $\hat{t}\in]-t,t[$. If $\vert \dddot{{\cal{E}}}\vert\leq c$, then for the critical domain 
\begin{eqnarray*}
{\cal{E}}(t)\geq {\cal{E}}(0)+\frac{t^2}{2}\left(\ddot{{\cal{E}}}(0)-c\:t\right)
\end{eqnarray*}
which shows that ${\cal{E}}(0)$ is minimal for small $t$.
In a first step we find upper bounds for $\dddot{m}, \dddot{J}$ and $\dddot{A}$. The main tool is the formula
\begin{eqnarray*}
\det\left(Id+tA)\right)=\sum_{k=0}^{\infty}\frac{1}{k!}\left(-\sum_{j=1}^{\infty}\frac{(-1)^j}{j}tr((t\:A)^j)\right)^k,
\end{eqnarray*}
where $tr(A)$ denotes the trace of the matrix $A$. In our case the matrix $A$ will depend on $t$ as well:
\begin{eqnarray*}
A:=\tilde{A}+\frac{t}{2}\tilde{B},\quad\hbox{where}\quad\tilde{A}=D_v, \:\tilde{B}=D_w.
\end{eqnarray*}
We now assume 
\begin{eqnarray}\label{thirdass}
\| D_v\|_{L^{\infty}(\Omega)}+\| D_w\|_{L^{\infty}(\Omega)}\leq 1\quad\hbox{and}\quad 0\leq t<1/2.
\end{eqnarray}
Then
\begin{eqnarray*}
\det\left(Id+tA)\right)\leq e^{\frac{c}{1-t}}\leq c,
\end{eqnarray*}
and $c$ does not depend on $v$ or $w$. With these assumptions we also get the estimates
\begin{eqnarray*}
\left\vert\frac{d}{dt}\det\left(Id+tA\right)\right\vert,\:
\left\vert \frac{d^2}{dt^2}\det\left(Id+tA\right)\right\vert,\:
\left\vert \frac{d^3}{dt^3}\det\left(Id+tA\right)\right\vert\leq c.
\end{eqnarray*}
Again $c$ does not depend on $v$ or $w$. 
With this we can easily prove the following lemma.
%%%%%%%%%%%%%%%%%%%%%%%%%%%%%%%%%%%%%%%
\begin{lemma}\label{dddotm}
Let $J$ (resp. $m$ and $A$) be defined as in \eqref{jform} (resp. \eqref{Wurzelk} and \eqref{aij1}).
Moreover, we assume \eqref{thirdass} for the vector fields $v$, $w$ and the parameter 
$t$.
Then the following estimates hold:
\begin{eqnarray*}
m(t)+\dot{m}(t)+\ddot{m}(t)+\dddot{m}(t)&\leq&c_0\\
J(t)+\dot{J}(t)+\ddot{J}(t)+\dddot{J}(t)&\leq&c_1\\
\|A(t)\|+\|\dot{A}(t)\|+\|\ddot{A}(t)\|+\|\dddot{A}(t)\|&\leq&c_1,
\end{eqnarray*}
where $c_0$ and $c_1$ do not depend on $v$ and $w$. $\|\cdot\|$ denotes any matrix norm.
\end{lemma}
%%%%%%%%%%%%%%%%%%%%%%%%%%%%%%%%%%%%%%%%%
In a final step we assume that for some number $c\in\rz$ we have
\begin{eqnarray}\label{thirduass}
\vert G'(u)\vert\leq c.
\end{eqnarray}
Then from \eqref{transsol1g} and \eqref{transsol2g} and the corresponding equation for $\dot{\tilde{u}}$ we get
\begin{eqnarray*}
&&\int_{\Omega}\vert\nabla\tilde{u}\vert^2\:dx+\alpha\oint_{\Omega}\tilde{u}^2\:dS\leq c\\
&&\int_{\Omega}\vert\nabla\dot{\tilde{u}}\vert^2\:dx+\alpha\oint_{\Omega}\dot{\tilde{u}}^2\:dS\leq c.
\end{eqnarray*}
%%%%%%%%%%%%%%%%%%%%%%%%%%%%%%%%%%%%%%%%%
\begin{theorem}
Let $v$ and $w$ be two smooth vector fields satisfying \eqref{thirdass}. Then there exists a number $c\in\rz$, which is independent of $v$ and $w$, such that
\begin{eqnarray*}
\vert\dddot{{\cal{E}}}(t)\vert\leq c\qquad\forall\:0\leq t<\frac{1}{2}.
\end{eqnarray*}
Consequently
\begin{eqnarray*}
\ddot{{\cal{E}}}(t)\geq \ddot{{\cal{E}}}(0) -c\:t \qquad\forall\:0\leq t<\frac{1}{2}.
\end{eqnarray*}
For $t$ sufficiently small we get the uniform positivity of $\ddot{{\cal{E}}}(t)$.
\end{theorem}
%%%%%%%%%%%%%%%%%%%%%%%%%%%%%%%%%%%%%
Since $\ddot{{\cal{E}}}(t)$ does not depend on $\ddot{\tilde{u}}$, it is also independent of the tangential component of $v$, and $w$.
%%%%%%%%%%%%%%%%%%%%%%%%%%%%%%%%%%%%%%%%%%
%%%%%%%%%%%%%%%%%%%%%%%%%%%%%%%%%%%%%%%%%%
\section{Back to Garabedian and Schiffer's second variation}\label{garab}
%%%%%%%%%%%%%%%%%%%%%%%%%%%%%%%%%%%%%%%%%
%%%%%%%%%%%%%%%%%%%%%%%%%%%%%%%%%%%%%%%%%
Garabedian and Schiffer computed the second domain variation of the first Dirichlet eigenvalue of the Laplace operator for the ball (\cite{GaSc53}). Since the Krahn - Faber inequality holds, one would expect the strict positivity of
$\ddot{\lambda}_{D}(0)$.
However, from the formula Garabedian and Schiffer obtained, namely
\begin{align*}
\ddot{\lambda}_D(0)= -\oint_{\p \Omega} (\p_\nu u)^2(v\cdot \nu)^2H\:dS- 2\int_{\Omega}\left(|\nabla \dot{\tilde{u}}(0)|^2-\lambda (  \dot{\tilde{u}}(0))^2\right)\:dx,
\end{align*}
it seems to be difficult to show that $\ddot{\lambda}_{D}(0)\geq0$.  Throughout this section we shall assume that $\int_{B_R}u^2\:dx=1$. Following the device of our paper we find
\begin{eqnarray*}
\frac{1}{2}\ddot{\lambda}_{D}(0)=\int_{B_R}\vert\nabla u'\vert^2-\lambda_{D}u'^2\:dx
+
\frac{n-1}{R}\oint_{\p B_R}u'^2\:dS.
\end{eqnarray*}
Here $u'$ satisfies the equation
\begin{eqnarray*}
\Delta u' +\lambda_{D}\:u'&=&0\qquad\hbox{in}\:\:B_R\\
u' &=&-(v\cdot\nu)\partial_{\nu}u\qquad\hbox{in}\:\:\partial B_R,
\end{eqnarray*}
where $u=u(|x|)$.
In this computation $\dot{\lambda}_{D}(0)=0$ is already taken into account.
We define
\begin{eqnarray*}
{\cal{R}}_{s}(\phi):=\frac{\int_{B_R}\vert\nabla \phi\vert^2-\lambda_{D}\phi^2\:dx}{\oint_{\p B_R}\phi^2\:dS},
\end{eqnarray*}
and we set
\begin{eqnarray*}
\mu:=\inf\left\{{\cal{R}}_{s}(\phi), \oint_{\p B_R} \phi\:dS=0\right\}.
\end{eqnarray*}
From the previous considerations we observe that $\mu=\mu_2-\alpha$, where $\mu_2$ is as in the previous subsection. As before
\begin{eqnarray*}
u=u(r)=c\:r^{\frac{2-n}{2}}\:J_{\frac{n-2}{2}}(\sqrt{\lambda_D}\:r).
\end{eqnarray*}
$\lambda_D$ is determined by the boundary condition $u(R)=0$, i.e.
\begin{eqnarray}\label{zero}
J_{\frac{n-2}{2}}(\sqrt{\lambda_D}R)=0.
\end{eqnarray}
By the same arguments as in the previous section
\begin{align*}
\frac{1}{2}\ddot{\lambda}_{D}(0)\geq\left\{ \frac{n}{R}-\frac{\sqrt{\lambda_D}J_{\frac{n}{2}+1}(\sqrt{\lambda_D}R)}{J_{\frac{n}{2}}(\sqrt{\lambda_D}R)}\right\}\oint_{\p B_R}u'^2\:dS.
\end{align*}
The identity \eqref{Bessel} and \eqref{zero} imply that
\begin{align*}
\ddot{\lambda}(0)\geq 0.
\end{align*}
As in the last section the equality sign can be excluded if $v$ satisfies \eqref{ortho1}.
%%%%%%%%%%%%%%%%%%%%%%%%%%%%%%%%%%%%%%%%%
\subsection{The Case of Dirichlet data}
\label{subdd}
%%%%%%%%%%%%%%%%%%%%%%%%%%%%%%%%%%%%%%%%%
In case of Dirichlet data $u=0$ on $\partial\Omega$ the energy ${\cal{E}}(t)$ has the form
\begin{eqnarray}\label{energyD}
{\cal{E}}(t)=\int_{\Omega}\nabla\tilde{u}A\nabla\tilde{u}\:dx-2\int_{\Omega}G(\tilde{u})\:J\:dx.
\end{eqnarray}
As in \eqref{transsol1g} the function $\tilde{u}$ solves
\begin{eqnarray*}
L_{A}\tilde{u}(t)+g(\tilde{u}(t))\:J(t) = 0\qquad\hbox{in}\quad\Omega
\end{eqnarray*}
with the boundary condition \eqref{transsol2g} replaced by
\begin{eqnarray}\label{transsolDg}
\tilde{u}(t) = 0\qquad\hbox{in}\quad\partial\Omega.
\end{eqnarray}
The function $u'$ solves
\begin{eqnarray}\label{uprimD}
\label{uprimD1}\Delta u' +g'(u)\:u'&=&0\qquad\hbox{in}\quad\Omega\\
\label{uprimD2}u'&=&-v\cdot\nabla u\qquad\hbox{in}\quad\partial\Omega.
\end{eqnarray}
In complete analogy with Chapter 4.1 we get
\begin{eqnarray}\label{edotD}
\dot{{\cal{E}}}(0)
=
\oint_{\partial\Omega}(v\cdot\nu)\left\{\vert\nabla u\vert^2-2G(u)\right\}\:dS.
\end{eqnarray}
Thus, any critical domain for which $\dot{{\cal{E}}}(0)=0$ satisfies the overdetermined boundary condition
\begin{eqnarray}\label{edotDeq}
\vert\nabla u\vert^2-2G(u)=const.\qquad\hbox{in}\quad\partial\Omega.
\end{eqnarray}
Hence
\begin{eqnarray}\label{overdetD}
\vert\nabla u\vert=c_0\qquad\hbox{in}\quad\partial\Omega. 
\end{eqnarray}
Note that, by a result of Serrin for positive solutions, this would already imply that 
$\Omega$ is a ball. For the second variation we observe that only ${\cal{S}}$ and 
\begin{eqnarray*}
{\cal{F}}_1+{\cal{F}}_4+{\cal{F}}_3+{\cal{F}}_6
\end{eqnarray*} 
contribute. Hence
\begin{eqnarray*}
\ddot{{\cal{E}}}(0)=\ddot{\cal{S}}(0)+{\cal{F}}_1(0)+{\cal{F}}_4(0)+{\cal{F}}_3(0)+{\cal{F}}_6(0).
\end{eqnarray*}
Computations very similar to those in Chapter 6 lead to the following lemma.
%%%%%%%%%%%%%%%%%%%%%%%%%%%%%%%%%%%%%%%%%
\begin{lemma}\label{edotDlem}
Let $\Omega$ be a smooth domain and let ${\cal{E}}(t)$ be as in \eqref{energyD}. Let $u$ be a solution of $\Delta u+g(u)=0$ in $\Omega$ and $u=0$ on $\partial\Omega$. Let $u'$ be a solution of \eqref{uprimD1} - \eqref{uprimD2}. For any critical domain $\Omega$ in the sense that $\dot{{\cal{E}}}(0)=0$ we have
\begin{eqnarray*}
\ddot{{\cal{E}}}(0)=2Q_g(u')+g(0)\oint_{\partial\Omega}u'^2\:\partial_{\nu}u\:dS
+
2(n-1)\:\oint_{\partial\Omega}u'^2\:H\:dS,
\end{eqnarray*}
where $H$ denotes the mean curvature of $\partial\Omega$.
\end{lemma}
%%%%%%%%%%%%%%%%%%%%%%%%%%%%%%%%%%%%%%%%%
%%%%%%%%%%%%%%%%%%%%%%%%%%%%%%%%%%%%%%%%%
\bigskip
%%%%%%%%%%%%%%%%%%%%%%%%%%%%%%%%%%%%%%%%%
%%%%%%%%%%%%%%%%%%%%%%%%%%%%%%%%%%%%%%%%%
{\bf{Acknowledgements}}
The authors would like to thank J. Arrieta, E. Harrell, M. Pierre and the referee for drawing our attention to some of the more recent literature on the subject.
%%%%%%%%%%%%%%%%%%%%%%%%%%%%%%%%%%%%%%%%%
%%%%%%%%%%%%%%%%%%%%%%%%%%%%%%%%%%%%%%%%%

\end{document}